\documentclass[12pt]{article}

\usepackage{amsmath,amssymb,overcite}
\usepackage[mathscr]{eucal}
\usepackage{theorem}

\renewcommand{\thesection}{\Roman{section}}

\makeatletter
\renewcommand{\@biblabel}[1]{%
{\textsuperscript{#1}}}
\makeatother

\newtheorem{lemma}{\it Lemma}
\newtheorem{proposition}{\it Proposition}
\newtheorem{corollary}{Corollary}
{\theorembodyfont{\rmfamily}
\newtheorem{remark}{\sc Remark}
\newtheorem{example}{\sc Example}}

\begin{document}

\vspace*{\fill}

\noindent{\sffamily\bfseries\Huge Differential calculi on some quantum
prehomogeneous vector spaces}

\begin{quotation}
\noindent{\sffamily\large S. Sinel'shchikov$^\mathsf{a)}$, A.
Stolin$^\mathsf{b)}$, and L. Vaksman$^\mathsf{a)}$}
{\it
\\ $^\mathsf{a)}$Mathematics Division, Institute for Low Temperature Physics
\& Engineering, 47 Lenin Ave, 61103 Kharkov, Ukraine
\\ $^\mathsf{b)}$Chalmers Tekniska H\"ogskola, Mathematik, 412 96,
G\"oteborg, Sweden}
\end{quotation}

\bigskip

This paper is devoted to study of differential calculi over quadratic
algebras, which arise in the theory of quantum bounded symmetric domains.
We prove that in the quantum case dimensions of the homogeneous components
of the graded vector spaces of $k$-forms are the same as in the classical
case. This result is well-known for quantum matrices.

The quadratic algebras, which we consider in the present paper, are
q-analogues of the polynomial algebras on prehomogeneous vector spaces of
commutative parabolic type. This enables us to prove that the de Rham
complex is isomorphic to the dual of a quantum analogue of the generalized
Bernstein-Gelfand-Gelfand resolution.

\vspace*{\fill}
\newpage

\section{\sffamily\bfseries\!\!\!\!\!\!. INTRODUCTION}

\indent

In his works Rubenthaler introduces and studies prehomogeneous vector
spaces of commutative parabolic type \cite{Rub} . Their quantum analogues
are crucial in the theory of quantum bounded symmetric domains. More
exactly, an irreducible bounded symmetric domain admits a standard
realisation as a unit ball in some finite dimensional complex Banach space,
which is a prehomogeneous vector space of commutative parabolic type.

The simplest example of such a vector space is a vector space of complex
$m\times n$ matrices considered as a space of linear maps from
$\mathbb{C}^m\to\mathbb{C}^n$. The algebra of polynomial functions on this
space is a $U\mathfrak{sl}_{m+n}$-module algebra by the following reason.
There exists a natural correspondence between linear maps from
$\mathbb{C}^n\to\mathbb{C}^m$ and an open subset in
$\mathbb{G}\mathrm{r}(n,n+m)$, the Grassman manifold of $n$-dimensional
subspaces in $\mathbb{C}^{n+m}$. Namely, to any linear map we assign its
graph.

Since the group $S(GL_m\times GL_n)$ acts naturally on the space of
$m\times n$-matrices, we see that the latter vector space is a
prehomogeneous vector space of commutative parabolic type.

Note, that the algebra $U\mathfrak{sl}_{m+n}$ is much bigger than
$U\mathfrak{s}(\mathfrak{gl}_m\times\mathfrak{gl}_n)$. Therefore the study
of the $U sl_{m+n}$-symmetry (sometimes it is called a hidden symmetry)
gives much more information than the study of the obvious $S(GL_n\times
GL_n)$-symmetry. In the middle of 90's it became clear that the quantum
matrices also have a hidden symmetry because it is well known that the
quantum analogue of the polynomial algebra on the space of $m\times n$
matrices is a $U_q\mathfrak{sl}_{m+n}$-module algebra \cite{SV0} .

Now let us consider the general case. Let $\mathfrak{g}$ be a simple
complex finite dimensional Lie algebra, and let $\mathfrak{t}$ be its
reductive subalgebra of the same rank. Let $P(T)$ be the parabolic subgroup
of $G$ with the reductive part $T$ (here the Lie groups $G$, $T$ correspond
to $\mathfrak{g}$, $\mathfrak{t}$). Clearly, the polynomial algebra on
$G/P(T)$ has a structure of $U\mathfrak{t}$-module algebra, and it is not
difficult to see that it has a structure of $U\mathfrak{g}$-module algebra.
A quantum analogue of this fact is also true what was shown in \cite{SV}.

In our paper we obtain a number of results about the
$U_q\mathfrak{g}$-covariant differential calculi on the polynomial algebra
on quantum prehomogeneous vector spaces constructed in \cite{SV}. Our goal
is to produce new series of quantum Harish-Chandra modules using the method
of cohomological induction. As it was shown in \cite{BaEast}, the
generalized Bernstein-Gelfand-Gelfand (BGG) resolutions are extremely
important in this theory. It is known that the generalized BGG resolution
of the trivial $U\mathfrak{g}$-module is dual to the De Rham resolution. We
obtain a quantum analogue of this result which is new even in the special
case of quantum matrices.

First we prove that the Hilbert series for the spaces of differential forms
on the considered prehomogeneous spaces for classical and related to them
quantum groups coincide. For case of matrices this was shown in
\cite{Maltsiniotis} . Next we prove that the polynomial algebras on quantum
prehomogeneous vector spaces of commutative type are quadratic algebras.
Our proof does not require separate consideration of serial and exceptional
simple Lie algebras (cf. \cite{Jak-Hermit, KMT}).

Similar ideas were used in the recent works of I. Heckenberger and S.
Kolb\cite{HeckenbergerKolb_last} . In the main text we will compare our
methods and results.

Finally, we would like to note that a hidden symmetry can be used in the
study of determinantal varieties, as it was shown in the recent papers of
Enright and his coauthors \cite{EnrHunz, Enr_except}.

\section{\sffamily\bfseries\!\!\!\!\!\!. THE
\boldmath$U_q\mathfrak{g}$-MODULE ALGEBRA
$\pmb{\mathbb{C}}[\mathfrak{p}^-]_q$}

\indent

Consider a simple complex Lie algebra $\mathfrak{g}$ with the Chevalley
generators $\{H_i,E_i,F_i\}_{i=1,2,\ldots,l}$. Let $\mathfrak{h}$,
$\mathfrak{n}^+$, $\mathfrak{n}^-$ be the Lie subalgebras of $\mathfrak{g}$
generated by $\{H_i\}_{i=1,2,\ldots,l}$, $\{E_i\}_{i=1,2,\ldots,l}$, and
$\{F_i\}_{i=1,2,\ldots,l}$, respectively.

Let $\{\alpha_j\}_{j=1,2,\ldots,l}$ be a system of simple roots and
$\mathbf{a}=(a_{ij})_{i,j=1,2,\ldots,l}$ be the Cartan matrix:
$a_{ij}=\alpha_j(H_i)$. The set of fundamental weights
$\{\overline{\omega}_j\}_{j=1,2,\ldots,l}$ forms a basis in
$\mathfrak{h}^*\cong\mathbb{C}^l$ and one has
$\alpha_j=\sum_{i=1}^la_{ij}\overline{\omega}_i$. There exists a unique
sequence of coprime positive integers $d_1,d_2,\ldots,d_l$ such that
$$d_ia_{ij}=d_ja_{ji},\qquad i,j=1,2,\ldots,l,$$
and the bilinear form in $\mathfrak{h}^*$ given by
$$(\alpha_i,\alpha_j)=d_ia_{ij},\qquad i,j=1,2,\ldots,l,$$
is positive definite. In this setting,
$(\overline{\omega}_i,\alpha_j)=d_i\delta_{ij}$.

Recall the standard notation $Q=\bigoplus_{i=1}^l\mathbb{Z}\alpha_i$ for
the root lattice, $P=\bigoplus_{i=1}^l\mathbb{Z}\overline{\omega}_i$ for
the weight lattice, and also $Q_+=\bigoplus_{i=1}^l\mathbb{Z}_+\alpha_i$,
$P_+=\bigoplus_{i=1}^l\mathbb{Z}_+\overline{\omega}_i$.

\medskip

The maximal root of $\mathfrak{g}$ is given by a linear combination of its
simple roots
\begin{equation}\label{decomposition} \sum_{j=1}^ln_j\alpha_j,\qquad
n_j\in\mathbb{Z}_+.
\end{equation}
Fix a simple root $\alpha_{l_0}$ used in the above sum with the coefficient
1: $n_{l_0}=1$. The tables in \cite{Bou4-6_fr} allow one to distinguish
easily all simple roots of this kind.

Let $H_0$ be a linear combination of $H_1,H_2,\ldots,H_l$ determined by
$$
\alpha_j(H_0)=
\begin{cases}
2, & j=l_0
\\ 0, & j\ne l_0.
\end{cases}
$$
Introduce the notation $\mathfrak{p}^-$, $\mathfrak{k}$, $\mathfrak{p}^+$
for the eigenspaces of $\operatorname{ad}_{H_0}$ which correspond to the
eigenvalues $-2,0,2$, respectively. Thus we get a decomposition
$$\mathfrak{g}=\mathfrak{p}^-\oplus\mathfrak{k}\oplus\mathfrak{p}^+,$$
where $\mathfrak{p}^\pm$ are commutative Lie algebras. They are called
prehomogeneous vector spaces of commutative parabolic type \cite{Rub} .

We are interested in a quantum analogue for the algebra of differential
forms with polynomial coefficients on the vector space $\mathfrak{p}^-$.

Let $s={\rm card}(P/Q)$. We assume that all algebras in this paper are
unital and they are defined over the field $\mathbb{C}(q^{\frac1s})$, the
field of rational functions in $q^{\frac1s}$.

Consider the Drinfeld-Jimbo quantum universal enveloping algebra
$U_q\mathfrak{g}$. It is determined by its generators $\{K_i^{\pm
1},E_i,F_i\}_{i=1,2,\ldots,l}$ and the well-known relations \cite[p.
52]{Jant}. It is a Hopf algebra with the comultiplication
$$
\triangle(K_i^{\pm 1})=K_i^{\pm 1}\otimes K_i^{\pm 1},\qquad
\triangle(E_i)=E_i\otimes 1+K_i\otimes E_i,\qquad
\triangle(F_i)=F_i\otimes K_i^{-1}+1\otimes F_i,
$$
the antipode $S$, and the counity $\varepsilon$.

In what follows all the $U_q\mathfrak{g}$-modules are left
$U_q\mathfrak{g}$-modules, unless the contrary is stated explicitly.

Let $q_j=q^{d_j}$. All the $U_q\mathfrak{g}$-modules, which we consider,
are weight modules:
$$
V=\bigoplus_{\lambda\in\mathbb{Z}^l}V_\lambda,\qquad V_\lambda=\left\{v\in
V\left|\:K_j^{\pm 1}v=q_j^{\pm\lambda_j}v,\;j=1,2,\ldots,l\right.\right\}.
$$

This allows one to introduce the linear operators
$\{H_j\}_{j=1,2,\ldots,l}$, (hence also $H_0$) as follows:
$$H_j|_{V_\lambda}=\lambda_j,\qquad\lambda\in\mathbb{Z}^l,$$
and to consider the following grading on $V$:
$$V=\bigoplus_rV[r],\qquad V[r]=\{v\in V|\:H_0v=2rv\}.$$

In what follows all the $U_q\mathfrak{g}$-modules are equipped with this
grading.

By definition we set $V^*=\bigoplus_rV[r]^*$ in the category of graded
vector spaces. $V^*$ has the following $U_q\mathfrak{g}$-module structure:
$$\xi l(v)=l(S(\xi)v),\qquad\xi\in U_q\mathfrak{g},\;v\in V,\;l\in V^*.$$

We need the generalized Verma modules $N(\mathfrak{q}^+,\lambda)$, which
are defined as follows. Let $U_q\mathfrak{q}^+$ be the Hopf subalgebra
generated by $\{K_i^{\pm 1},E_i,F_i\}_{i\ne l_0}\cup\{K_{l_0}^{\pm
1},E_{l_0}\}$ and
$$
\mathscr{P}_+=\mathbb{Z}_+^{l_0-1}\times\mathbb{Z}\times\mathbb{Z}_+^{l-l_0}
\hookrightarrow
P,\qquad(n_1,n_2,\ldots,n_l)\mapsto\sum_{j=1}^ln_j\overline{\omega}_j,
$$
$\lambda\in\mathscr{P}_+$. Consider the simple finite dimensional
$U_q\mathfrak{q}^+$-module $L(\mathfrak{q}^+,\lambda)$ with the highest
weight $\lambda$, together with the induced $U_q\mathfrak{g}$-module
$N(\mathfrak{q}^+,\lambda)=U_q\mathfrak{g}\otimes_{U_q\mathfrak{q}^+}
L(\mathfrak{q}^+,\lambda)$. It is easy to show that
$N(\mathfrak{q}^+,\lambda)$ admits a description in terms of a generator
$v(\mathfrak{q}^+,\lambda)$ and the defining relations
\begin{gather*}
E_iv(\mathfrak{q}^+,\lambda)=0,\qquad K_i^{\pm
1}v(\mathfrak{q}^+,\lambda)=q_i^{\pm\lambda_i}v(\mathfrak{q}^+,\lambda),
\qquad i=1,2,\ldots,l,
\\ F_j^{\lambda_j+1}v(\mathfrak{q}^+,\lambda)=0,\qquad j\ne l_0;
\end{gather*}
for a similar result in the classical case $q=1$ see \cite{Lep} .

Turn to producing a $U_q\mathfrak{g}$-module algebra
$\mathbb{C}[\mathfrak{p}^-]_q$ that will work as a $q$-analogue of the
$U\mathfrak{g}$-module algebra of polynomials on the vector space
$\mathfrak{p}^-$. Consider the graded vector space $N(\mathfrak{q}^+,0)$
together with the dual graded vector space
$\mathbb{C}[\mathfrak{p}^-]_q=N(\mathfrak{q}^+,0)^*$. We equip
$N(\mathfrak{q}^+,0)$ with a structure of graded coalgebra and
$\mathbb{C}[\mathfrak{p}^-]_q$ with a structure of graded algebra by
duality.

This approach was used by Drinfeld when he constructed algebras of
functions on quantum groups. It is based on the fact that the vector space
that is dual to a coalgebra is an algebra. Note that the tensor factors
when passing to dual spaces remain non-permuted. On the contrary, in the
tensor category of $U_q\mathfrak{g}$-modules one has $V_1^*\otimes
V_2^*\hookrightarrow(V_2\otimes V_1)^*$, that is, the tensor factors change
their places under passing to dual spaces. This inconsistence can be
avoided by using the Hopf algebra $U_q\mathfrak{g}^\mathrm{cop}$ which is
just $U_q\mathfrak{g}$, but its comultiplication is replaced with the
opposite one.

Let us treat the generalized Verma modules $N(\mathfrak{q}^+,\lambda)$ as
the objects in the tensor category of
$U_q\mathfrak{g}^\mathrm{cop}$-modules, and the corresponding dual graded
modules over the quantum universal enveloping algebra as the objects in the
tensor category of $U_q\mathfrak{g}$-modules. This convention yields a
canonical isomorphism
$$
(N(\mathfrak{q}^+,\lambda)\otimes N(\mathfrak{q}^+,\mu))^*\cong
N(\mathfrak{q}^+,\lambda)^*\otimes N(\mathfrak{q}^+,\lambda)^*,
$$
which allows one to avoid permutation of the tensor factors.

The morphism
$$
\triangle_0:N(\mathfrak{q}^+,0)\to N(\mathfrak{q}^+,0)\otimes
N(\mathfrak{q}^+,0),\qquad\triangle_0:v(\mathfrak{q}^+,0)\mapsto
v(\mathfrak{q}^+,0)\otimes v(\mathfrak{q}^+,0)
$$
in the tensor category of $U_q\mathfrak{g}^\mathrm{cop}$-modules equips the
generalized Verma module $N(\mathfrak{q}^+,0)$ with a structure of graded
$U_q\mathfrak{g}^\mathrm{cop}$-module coalgebra.

The adjoint linear map
$m:\mathbb{C}[\mathfrak{p}^-]_q\otimes\mathbb{C}[\mathfrak{p}^-]_q\to
\mathbb{C}[\mathfrak{p}^-]_q$ to $\triangle_0$ is a morphism of
$U_q\mathfrak{g}$-modules and equips $\mathbb{C}[\mathfrak{p}^-]_q$ with a
structure of graded $U_q\mathfrak{g}$-module algebra (see \cite{SV}~).

So we obtain a q-analogue for the polynomial algebra on $\mathfrak{p}^-$.

One should note that in the papers by Joseph and his coauthors (see, e.g.,
\cite{Kebe}) a more general class of $U_q\mathfrak{g}$-module algebras is
treated.

\section{\sffamily\bfseries\!\!\!\!\!\!. COVARIANT DIFFERENTIAL CALCULI}
\label{calculi}

\indent

In the classical case $q=1$, the linear map that is adjoint to a morphism
of the generalized Verma modules appears to be a covariant differential
operator on $\mathfrak{p}^-$, see \cite{Har-Jakob} . This fact can be used
as a hint to apply a duality argument in the quantum case to produce a
first order differential calculus over $\mathbb{C}[\mathfrak{p}^-]_q$, see
\cite{SV} . Let us recall some definitions \cite{KlSch} .

Let $F$ be an algebra. A first order differential calculus over $F$ is an
$F$-bimodule $M$ together with a linear map $d:F\to M$ such that:
\begin{enumerate}
\item for all $f_1,f_2\in F$
\begin{equation}\label{4.1.1}
d(f_1\cdot f_2)=df_1\cdot f_2+f_1\cdot df_2;
\end{equation}
\item $M$ is a linear span of the vectors $f_1\cdot df_2\cdot f_3$ with
$f_1,f_2,f_3\in F$.
\end{enumerate}

Let $A$ be a Hopf algebra and let $F$ be an $A$-module algebra. A first
order differential calculus $(M,d)$ over $F$ is called covariant if $M$ is
an $A$-module $F$-bimodule and $d$ is a morphism of $A$-modules. An
isomorphism of such calculi is introduced in a natural way.

Our next step is to produce a covariant first order differential calculus
over $\mathbb{C}[\mathfrak{p}^-]_q$.

Given any $\lambda\in\mathscr{P}_+$, the generalized Verma module
$N(\mathfrak{q}^+,\lambda)$ is a $U_q\mathfrak{g}^\mathrm{cop}$-module
$N(\mathfrak{q}^+,0)$-bicomodule:
\begin{gather*}
N(\mathfrak{q}^+,\lambda)\to N(\mathfrak{q}^+,0)\otimes
N(\mathfrak{q}^+,\lambda),\qquad v(\mathfrak{q}^+,\lambda)\mapsto
v(\mathfrak{q}^+,0)\otimes v(\mathfrak{q}^+,\lambda),
\\ N(\mathfrak{q}^+,\lambda)\to N(\mathfrak{q}^+,\lambda)\otimes
N(\mathfrak{q}^+,0),\qquad v(\mathfrak{q}^+,\lambda)\mapsto
v(\mathfrak{q}^+,\lambda)\otimes v(\mathfrak{q}^+,0).
\end{gather*}
Hence the dual graded vector space is a $U_q\mathfrak{g}$-module
$\mathbb{C}[\mathfrak{p}^-]_q$-bimodule. In particular,
$$
\Lambda^1(\mathfrak{p}^-)_q\overset{\mathrm{def}}{=}
N(\mathfrak{q}^+,-\alpha_{l_0})^*
$$
is a $U_q\mathfrak{g}$-module $\mathbb{C}[\mathfrak{p}^-]_q$-module. It is
a q-analogue for the space of 1-forms with polynomial coefficients.

Here the linear operator that is dual to the following morphism of the
generalized Verma modules
$$
N(\mathfrak{q}^+,-\alpha_{l_0})\to N(\mathfrak{q}^+,0),\qquad
v(\mathfrak{q}^+,-\alpha_{l_0})\mapsto F_{l_0}v(\mathfrak{q}^+,0),
$$
is called a differential. Our definitions imply \eqref{4.1.1}, together
with the claim that the differential
$d:\mathbb{C}[\mathfrak{p}^-]_q\to\Lambda^1(\mathfrak{p}^-)_q$ is a
morphism of $U_q\mathfrak{g}$-modules. In the sequel it is shown that
$\Lambda^1(\mathfrak{p}^-)_q$ is a linear span of $\{f_1\cdot df_2\cdot
f_3|\:f_1,f_2,f_3\in\mathbb{C}[\mathfrak{p}^-]_q\}$.

\medskip

The results obtained by Heckenberger and Kolb \cite{HeckenbergerKolb04}
allow one to assume that, under some reasonable irreducibility assumptions,
this covariant first order differential calculus is unique up to
isomorphism.

\medskip

Consider a graded algebra
$\Omega=\bigoplus\limits_{i\in\mathbb{Z}_+}\Omega_i$ along with a linear
degree 1 map $d:\Omega\to\Omega$. The pair $(\Omega,d)$ is called a
differential graded algebra if $d^2=0$ and
\begin{equation}\label{4.1.15}
d(\omega'\cdot\omega'')=d\omega'\cdot\omega''+(-1)^n\omega'\cdot
d\omega'',\qquad\omega'\in\Omega_n,\,\omega''\in\Omega.
\end{equation}

A differential calculus over an algebra $F$ is a differential graded
algebra $(\Omega,d)$ such that $\Omega_0=F$, and $\Omega$ is generated by
$\Omega_0\oplus d\Omega_0$.

Suppose that a first order differential calculus $(M,d)$ over $F$ is given.
By definition, the associated universal differential calculus
$(\Omega^\mathrm{univ},d^\mathrm{univ})$ over this algebra should have the
following properties:
\begin{enumerate}
\item $\Omega_1^\mathrm{univ}=M$; \item $d^\mathrm{univ}|_F=d$;

\item given a differential calculus $(\Omega',d')$ over $F$ which
    satisfies the above two properties ($\Omega_1'=M$, $d'|_F=d$),
    there exists a homomorphism of differential graded algebras
    $\Omega^\mathrm{univ}\to\Omega'$ which is the identity when
    restricted to $\Omega_0^\mathrm{univ}\oplus\Omega_1^\mathrm{univ}$.
\end{enumerate}

Consider a Hopf algebra $A$ and an $A$-module algebra $F$. A differential
calculus $(\Omega,d)$ is said to be covariant if $\Omega$ is an $A$-module
algebra and $d$ is an endomorphism of the $A$-module $\Omega$. It is known
\cite{KlSch} , p. 463 -- 464, that the universal differential calculus
exists and is unique up to isomorphism. It is covariant if the original
first order differential calculus possesses this property.

\medskip

{\sc Remark.} The covariance notion we use here is more general than that
in \cite{KlSch} where it is implicit that all the considered
$U_q\mathfrak{g}$-modules are $U_q\mathfrak{g}$-finite. This difference
does not affect the proof of covariance for the universal differential
calculus expounded in \cite{KlSch} , p. 464.

In what follows, we obtain a series of results on the universal
differential calculus $\left(\Lambda(\mathfrak{p}^-)_q,d\right)$ of the
first order differential calculus
$\left(\Lambda^1(\mathfrak{p}^-)_q,d\right)$.

\section{\sffamily\bfseries\!\!\!\!\!\!. PBW-BASES AND R-MATRICES}
\label{PBW-R}

\indent

Recall the standard notation and some well-known results. Let
$U_q\mathfrak{h}$, $U_q\mathfrak{n}^+$, $U_q\mathfrak{n}^-$ be the
subalgebras of $U_q\mathfrak{g}$ generated by $\{K_i^{\pm
1}\}_{i=1,2,\ldots,l}$, $\{E_i\}_{i=1,2,\ldots,l}$, and
$\{F_i\}_{i=1,2,\ldots,l}$, respectively. The linear maps
\begin{align*}
U_q\mathfrak{n}^-\otimes U_q\mathfrak{h}\otimes U_q\mathfrak{n}^+ &\to
U_q\mathfrak{g}, & u^-\otimes u^0\otimes u^+ &\mapsto u^-u^0u^+, &
\\ U_q\mathfrak{n}^+\otimes U_q\mathfrak{h}\otimes U_q\mathfrak{n}^- &\to
U_q\mathfrak{g}, & u^+\otimes u^0\otimes u^- &\mapsto u^+u^0u^- &
\end{align*}
are the isomorphisms of vector spaces \cite{Jant} , p. 66.

We are to construct bases for the vector spaces $U_q\mathfrak{h}$,
$U_q\mathfrak{n}^\pm$, together with the associated bases for
$U_q\mathfrak{g}$ that are similar to the Poincar\'e-Birkhoff-Witt bases.
Obviously, the elements $K_1^{j_1}K_2^{j_2}\dots K_l^{j_l}$, with
$j_1,j_2,\ldots,j_l\in\mathbb{Z}$, form a basis of $U_q\mathfrak{h}$. What
remains is to produce the bases of $U_q\mathfrak{n}^+$,
$U_q\mathfrak{n}^-$.

Obviously we have
$$
\mathfrak{h}^*\overset{\cong}{\longrightarrow}\mathbb{C}^l,\qquad
\lambda\mapsto(\lambda(H_1),\lambda(H_2),\ldots,\lambda(H_l)).
$$

The decomposition
$\mathfrak{g}=\mathfrak{n}^-\oplus\mathfrak{h}\oplus\mathfrak{n}^+$
determines an associated decomposition $\Phi=\Phi^+\cup\Phi^-$ for the set
of roots of $\mathfrak{g}$.

The Weyl group $W$ is generated by the simple reflections
$$s_i:\lambda\mapsto\lambda-\lambda(H_i)\alpha_i,\qquad i=1,2,\ldots,l,$$
in $\mathfrak{h}^*$. Choose a reduced expression
\begin{equation}\label{A.2.14}
w_0=s_{i_1}s_{i_2}s_{i_3}\ldots s_{i_M}
\end{equation}
of the longest element $w_0$ of $W$. It determines an associated linear
order
\begin{equation}\label{A.2.15}
\beta_1=\alpha_{i_1},\qquad\beta_2=s_{i_1}(\alpha_{i_2}),\qquad\ldots,\qquad
\beta_M=s_{i_1}s_{i_2}\ldots s_{i_{M-1}}(\alpha_{i_M})
\end{equation}
on $\Phi^+$. Let us use the Lusztig braid group representation by the
automorphisms $T_1,T_2,\ldots,T_n$ of the algebra $U_q\mathfrak{g}$ to
obtain ``q-analogues of the root vectors'' $E_{\beta_k}$, $F_{\beta_k}$
(see, e.g., \cite{Rosso, Jant, DeConProc}) . The following proposition was
proved by Lusztig.

\begin{proposition}
The monomials
\begin{equation}\label{A.2.17}
E_{\beta_1}^{j_1}E_{\beta_2}^{j_2}\ldots E_{\beta_M}^{j_M},\qquad
j_1,j_2,\ldots j_M\in\mathbb{Z}_+,
\end{equation}
form a basis in $U_q\mathfrak{n}^+$, and the monomials
\begin{equation}\label{A.2.18}
F_{\beta_M}^{j_M}F_{\beta_{M-1}}^{j_{M-1}}\ldots F_{\beta_1}^{j_1},\qquad
j_1,j_2,\ldots j_M\in\mathbb{Z}_+,
\end{equation}
form a basis in $U_q\mathfrak{n}^-$.
\end{proposition}

Further we present here the commutation relations obtained by Levendorski
and Soibelman \cite{DeConProc} , p. 67--68.

\begin{proposition}\label{2.2.2}
\begin{enumerate}
\item \ \ For all $i<j$
\begin{eqnarray}
E_{\beta_i}E_{\beta_j}-q^{(\beta_i,\beta_j)}E_{\beta_j}E_{\beta_i} &=&
\sum_{\mathbf{m}\in\mathbb{Z}_+^M}C_\mathbf{m}'(q)\cdot
E^\mathbf{m},\label{A.2.19}
\\ F_{\beta_i}F_{\beta_j}-q^{-(\beta_i,\beta_j)}F_{\beta_j}F_{\beta_i} &=&
\sum_{\mathbf{m}\in\mathbb{Z}_+^M}C_\mathbf{m}''(q)\cdot
F^\mathbf{m},\label{A.2.20}
\end{eqnarray}
with $\mathbf{m}=(m_1,m_2,\ldots,m_M)$,
$E^\mathbf{m}=E_{\beta_1}^{m_1}E_{\beta_2}^{m_2}\ldots
E_{\beta_M}^{m_M}$,
$F^\mathbf{m}=F_{\beta_M}^{m_M}F_{\beta_{M-1}}^{m_{M-1}}\ldots
F_{\beta_1}^{m_1}$. The coefficients $C_\mathbf{m}'(q)$,
$C_\mathbf{m}''(q)$ can be non-zeros only when $m_1=m_2=\ldots=m_i=0$
and $m_j=m_{j+1}=\ldots=m_M=0$.

\item \ \ $C_\mathbf{m}'(q),\, C_\mathbf{m}''(q)\in \mathbb{Q}[q,q^{-1}]$.
\end{enumerate}
\end{proposition}

\begin{corollary}
The algebra $U_q\mathfrak{g}$ is a domain.
\end{corollary}

\begin{corollary}
The monomials
\begin{equation*}
E_{\beta_M}^{j_M}E_{\beta_{M-1}}^{j_{M-1}}\ldots E_{\beta_1}^{j_1},\qquad
j_1,j_2,\ldots j_M\in\mathbb{Z}_+,
\end{equation*}
form a basis in $U_q\mathfrak{n}^+$, and the monomials
\begin{equation*}
F_{\beta_1}^{j_1}F_{\beta_2}^{j_2}\ldots F_{\beta_M}^{j_M},\qquad
j_1,j_2,\ldots j_M\in\mathbb{Z}_+,
\end{equation*}
form a basis in $U_q\mathfrak{n}^-$.
\end{corollary}

Consider the Hopf algebras $U_q\mathfrak{b}^+$, $U_q\mathfrak{b}^-$
generated by $\{K_i^{\pm 1},E_i\}_{i=1,2,\ldots,l}$ and $\{K_i^{\pm
1},F_i\}_{i=1,2,\ldots,l}$, respectively. The $U_q\mathfrak{g}$-modules
discussed in this paper are weight and $U_q\mathfrak{b}^+$-finite (or
$U_q\mathfrak{b}^-$-finite), that is, $\dim(U_q\mathfrak{b}^+\cdot
v)<\infty$ for each vector $v$ (respectively, $\dim(U_q\mathfrak{b}^-\cdot
v)<\infty$).

Consider a pair of weight $U_q\mathfrak{g}$-modules $V'$, $V''$. Assume
that either $V'$ is $U_q\mathfrak{b}^+$-finite or $V''$ is
$U_q\mathfrak{b}^-$-finite. In this context Drinfeld introduced the
operators
$$\check{R}_{V',V''}:V'\otimes V''\to V''\otimes V',$$
which are quantum analogues of the naive permutation
$$
\sigma_{V',V''}:V'\otimes V''\to V''\otimes
V',\qquad\sigma_{V',V''}:v'\otimes v''\mapsto v''\otimes v'
$$
of tensor factors \cite{ChP, KlSch} . Here are some properties of the maps
$\check{R}_{V',V''}$, which will be used in the sequel.

\begin{proposition}\label{A.2.8}
Let $V'$, $V''$, $W'$, $W''$ be weight $U_q\mathfrak{b}^+$-finite
$U_q\mathfrak{g}$-modules (or $U_q\mathfrak{b}^-$-finite
$U_q\mathfrak{g}$-modules) and $f':V'\to W'$, $f'':V''\to W''$ be the
morphisms of $U_q\mathfrak{g}$-modules.
\begin{enumerate}
\item The linear map $\check{R}_{V',V''}$ is invertible and is a
    morphism of $U_q\mathfrak{g}$-modules.

\item $(f''\otimes
f')\cdot\check{R}_{V',V''}=\check{R}_{W',W''}\cdot(f'\otimes f'')$.
\end{enumerate}
\end{proposition}

\begin{proposition}\label{A.2.9}
Let $V$, $V'$, $V''$ be weight $U_q\mathfrak{b}^+$-finite
$U_q\mathfrak{g}$-modules (or $U_q\mathfrak{b}^-$-finite
$U_q\mathfrak{g}$-modules).
\begin{enumerate}
\item $\check{R}_{V'\otimes V'',V}=
\left(\check{R}_{V',V}\otimes\operatorname{id}_{V''}\right)
\left(\operatorname{id}_{V'}\otimes\check{R}_{V'',V}\right)$,

\item $\check{R}_{V,V'\otimes V''}=
\left(\operatorname{id}_{V'}\otimes\check{R}_{V,V''}\right)
\left(\check{R}_{V,V'}\otimes\operatorname{id}_{V''}\right)$,

\item
$\check{R}_{V,\mathbb{C}}=\check{R}_{\mathbb{C},V}=\operatorname{id}_{V}$.
\end{enumerate}
\end{proposition}

Recall an explicit form for $\check{R}_{V',V''}$. We intend to use
q-analogues for the root vectors $E_{\beta_i}$, $F_{\beta_i}$,
$i=1,2,\ldots,M$. Let
$$
\exp_q(t)=\sum_{i=0}^\infty\left(\prod_{j=1}^i\frac{1-q}{1-q^j}\right)t^i,
$$
$q_\beta=q^{\frac{(\beta,\beta)}2}$, and
$t_0\in\mathfrak{h}\otimes\mathfrak{h}$ be given by
$$(\lambda,\mu)=\lambda\otimes\mu(t_0),\qquad\lambda,\mu\in\mathfrak{h}^*.$$
Let $V'$, $V''$ be weight $U_q\mathfrak{b}^+$-modules and either $V'$ is
$U_q\mathfrak{b}_+$-finite or $V''$ is $U_q\mathfrak{b}^-$-finite. It is a
consequence of the results of \cite{KR, LS} that
\begin{equation}\label{R-check}
\check{R}_{V',V''}=\sigma_{V',V''}\cdot R_{V',V''},
\end{equation}
where $R_{V',V''}:V'\otimes V''\to V'\otimes V''$,
\begin{equation}\label{A.2.21}
R_{V',V''}\,(v'\otimes
v'')=\prod_{\beta\in\Phi^+}^\curvearrowleft\exp_{q_\beta^2}
\left(\left(q_\beta^{-1}-q_\beta\right)E_\beta\otimes
F_\beta\right)q^{-t_0}\,(v'\otimes v''),
\end{equation}
and the sign $\curvearrowleft$ indicates that the multipliers are written
in the decreasing order of indices
$$
\exp_{q_{\beta_M}^2}\left(\left(q_{\beta_M}^{-1}-q_{\beta_M}\right)
E_{\beta_M}\otimes F_{\beta_M}\right)\ldots
\exp_{q_{\beta_1}^2}\left(\left(q_{\beta_1}^{-1}-q_{\beta_1}\right)
E_{\beta_1}\otimes F_{\beta_1}\right)q^{-t_0}.
$$

Notice that
\begin{equation}\label{A.2.24}
t_0=\sum_k\frac{I_k\otimes I_k}{(I_k,I_k)}
\end{equation}
for any orthogonal basis $\{I_j\}_{j=1,2,\ldots,l}$ of
$\mathfrak{h}\cong\mathfrak{h}^*$.

It follows from Proposition 4.2 of \cite{Drinf2_engl} that
$$
\left(\check{R}_{N(\mathfrak{q}^+,\lambda)^*,
N(\mathfrak{q}^+,\mu)^*}\right)^*=
\left(\check{R}_{N(\mathfrak{q}^+,\lambda),
N(\mathfrak{q}^+,\mu)}\right)^{-1}.
$$
Here $N(\mathfrak{q}^+,\lambda)$, $N(\mathfrak{q}^+,\mu)$ are the objects
in the category of $U_q\mathfrak{g}^\mathrm{cop}$-modules, and
$N(\mathfrak{q}^+,\lambda)^*$, $N(\mathfrak{q}^+,\mu)^*$ are the graded
dual $U_q\mathfrak{g}$-modules.

\medskip

Turn to the finite dimensional weight $U_q\mathfrak{g}$-modules (these are
called type 1 $U_q\mathfrak{g}$-modules).

Let $\lambda\in P$. Just as in the classical case $q=1$, the Verma module
$M(\lambda)$ admits a description in terms of its generator $v(\lambda)$
and the defining relations
\begin{equation}\label{Verma-relations}
E_iv(\lambda)=0,\qquad K_i^{\pm
1}v(\lambda)=q_i^{\pm\lambda_i}v(\lambda),\qquad i=1,2,\ldots,l.
\end{equation}

The weight vectors
\begin{equation}\label{A.2.32}
v_J(\lambda)=F_{\beta_M}^{j_M}F_{\beta_{M-1}}^{j_{M-1}}\ldots
F_{\beta_1}^{j_1}v(\lambda),\qquad J=(j_1,j_2,\ldots,j_M)\in\mathbb{Z}_+^M,
\end{equation}
form a basis of the vector space $M(\lambda)$. Hence $M(\lambda)$ is a
weight $U_q\mathfrak{g}$-module with the same dimensions of weight
subspaces as in the classical case $q=1$.

The Verma module $M(\lambda)$ possesses the largest proper submodule
$K(\lambda)$. It is obvious that the quotient module
$L(\lambda)=M(\lambda)/K(\lambda)$ is simple.

$L(\lambda)$ is finite dimensional if and only if $\lambda\in P_+$. In this
case $K(\lambda)$ is the only proper submodule of finite codimension, and
$L(\lambda)$ admits a description in terms of its generator $v(\lambda)$,
the relations \eqref{Verma-relations}, together with the additional
relations $F_i^{\lambda_i+1}v(\lambda)=0$, $i=1,2,\ldots,l$.

The simple weight $U_q\mathfrak{g}$-modules $L(\lambda)$, $\lambda\in P_+$,
are pairwise non-isomorphic, and every simple weight finite dimensional
$U_q\mathfrak{g}$-module is isomorphic to one of them.

\begin{proposition}[cf.\cite{Jant} , p. 76 -- 77]
Given any non-zero element $\xi\in U_q\mathfrak{g}$, there exists
$\lambda\in P_+\cap Q$ such that $\xi L(\lambda)\ne 0$.
\end{proposition}

\begin{proposition}[\cite{Jant} , p. 81]\label{A.2.30}
The family of weights for $L(\lambda)$ and their multiplicities remain
intact under passage from $U\mathfrak{g}$ to $U_q\mathfrak{g}$, that is
from the classical case to the quantum case.
\end{proposition}

\begin{proposition}\label{100107}[\cite{Jant} , p. 82]
Every weight finite dimensional $U_q\mathfrak{g}$-module is semisimple.
\end{proposition}

It follows from the above results that
$$
L(\lambda)\otimes L(\mu)\approx\sum_{\nu\in P_+}c_{\lambda\mu}^\nu
L(\nu)\qquad\text{for all\ }\lambda,\mu\in P_+,
$$
and the multiplicities $c_{\lambda\mu}^\nu$ of $L(\nu)$ in
$L(\lambda)\otimes L(\mu)$ are the same as those in the classical case
$q=1$.

Let $P_{\lambda\mu}^\nu$ be the projection in $L(\lambda)\otimes L(\mu)$
onto the isotypic component which is multiple to $L(\nu)$ and parallel to
the sum of all other isotypic components.

\begin{proposition}[\cite{Drinf2_engl} , p. 333]\label{A.2.33}
For all $\lambda,\mu\in P_+$
\begin{equation}\label{A.2.43}
\check{R}_{L(\mu),L(\lambda)}\check{R}_{L(\lambda),L(\mu)}=\bigoplus_{\nu\in
P_+}q^{(\lambda,\lambda+2\rho)+(\mu,\mu+2\rho)-(\nu,\nu+2\rho)}
P_{\lambda\mu}^\nu.
\end{equation}
\end{proposition}

Below we sketch a standard method of reducing some problems related to
$U_q\mathfrak{g}$-modules, to the problems of the classical theory of
$U\mathfrak{g}$-modules.

The principal observation here is that many properties of
$U_q\mathfrak{g}$-modules can be formulated and proved in terms of their
distinguished submodules over the ring
$\mathscr{A}=\mathbb{Q}[q^{1/s},q^{-1/s}]$ of Laurent polynomials with
rational coefficients in the indeterminate $q^{1/s}$. In fact, consider the
$\mathscr{A}$-subalgebra $U_\mathscr{A}$ in $U_q\mathfrak{g}$ generated by
$K_i^{\pm 1}$, $E_i$, $F_i$, $h_i=\dfrac{K_i-K_i^{-1}}{q_i-q_i^{-1}}$,
$i=1,2,\ldots,l$. It is known from \cite{DeConKac} that a list of defining
relations between these generators can be derived from the standard list of
relations for $U_q\mathfrak{g}$ by replacing therein the relation
$E_iF_j-F_jE_i=\delta_{ij}\dfrac{K_i-K_i^{-1}}{q_i-q_i^{-1}}$ with
$E_iF_j-F_jE_i=\delta_{ij}h_i$ and adding the relation
$$(q_i-q_i^{-1})h_i=K_i-K_i^{-1},\qquad i=1,2,\ldots,l.$$
Certainly, the $\mathscr{A}$-subalgebra $U_\mathscr{A}$ inherits a
structure of the Hopf algebra.

Consider the homomorphisms $j:U_\mathscr{A}\to U\mathfrak{g}$,
$i:U_\mathscr{A}\to U_q\mathfrak{g}$ given by
\begin{gather*}
j(K_i^{\pm 1})=1,\quad j(E_i)=E_i,\quad j(F_i)=F_i,\quad j(h_i)=H_i,\qquad
i=1,2,\ldots,l,
\\ i(K_i^{\pm 1})=K_i^{\pm 1},\quad i(E_i)=E_i,\quad i(F_i)=F_i,\quad
i(h_i)=\dfrac{K_i-K_i^{-1}}{q_i-q_i^{-1}},\qquad i=1,2,\ldots,l.
\end{gather*}

The first homomorphism allows one to elaborate the classical theory of
$U\mathfrak{g}$-modules in studying $U_\mathscr{A}$-modules, and the second
one makes it possible to transfer the results related to
$U_\mathscr{A}$-modules onto $U_q\mathfrak{g}$-modules.

\section{\sffamily\bfseries\!\!\!\!\!\!. GENERATORS AND DEFINING RELATIONS
FOR THE ALGEBRA $\pmb{\mathbb{C}[\mathfrak{p}^-]_q}$}\label{new_item}

\indent

The main results of this section were obtained by Heckenberger and Kolb in
\cite{HeckenbergerKolb04} , while the auxiliary statements presented here
are new and will be used in the sequel.

We start with the well-known properties of the Weyl group $W$. Let
$\mathbb{S}=\{1,2,\ldots,l\}\setminus\{l_0\}$, $W_\mathbb{S}\subset W$ be
the subgroup generated by the simple reflections $s_i$, $i\in\mathbb{S}$,
and
$$
W^\mathbb{S}=\{w\in W\,|\,l(vw)\ge l(w)\;\text{for all}\;v\in
W_\mathbb{S}\,\}.
$$
It was shown by Kostant \cite{Kostant_1961} that any element $w\in W$ can
be represented uniquely as a product $w=w_\mathbb{S}\cdot w^\mathbb{S}$,
where $w_\mathbb{S}\in W_\mathbb{S}$, $w^\mathbb{S}\in\,W^\mathbb{S}$, and
$l(w)=l(w_\mathbb{S})+l(w^\mathbb{S})$. In particular, one has
$$
w_0\,=\,w_{0,\mathbb{S}}\cdot\;w_0^\mathbb{S},\qquad
l(w_0)=l(w_{0,\mathbb{S}})\,+\,l(w_0^\mathbb{S})
$$
for the longest element $w_0$ of the Weyl group $W$. In the above setting,
$w_{0,\mathbb{S}}$ is the longest element of the subgroup $W_\mathbb{S}$.

Fix the reduced expressions
\begin{equation}\label{2.1.3}
w_{0,\mathbb{S}}\,=\,s_{i_1}s_{i_2}\ldots s_{i_{M'}},\qquad
w_0^\mathbb{S}\,=\,s_{i_{M'+1}}s_{i_{M'+2}}\ldots s_{i_M}.
\end{equation}
Their concatenation $w_0=s_{i_1}s_{i_2}\ldots s_{i_M}$ is a reduced
expression for $w_0$.

Use it to produce a basis of $U_q\mathfrak{g}$. The algebra
$U_q\mathfrak{g}$ is a free right $U_q\mathfrak{q}^+$-module with the basis
$$
F_{\beta_M}^{j_M}F_{\beta_{M-1}}^{j_{M-1}}\ldots
F_{\beta_{M'+1}}^{j_{M'+1}},\qquad(j_M,j_{M-1},\ldots,j_{M'+1})\in
\mathbb{Z}_+^{M-M'}.
$$
Thus one has

\begin{proposition}\label{2.1.1}
Let $\lambda\in\mathscr{P}_+$ and $\{v_1,v_2,\ldots,v_d\}$ be a basis of
the vector space
$L(\mathfrak{q}^+,\lambda)$. Then the homogeneous elements
\begin{equation}\label{2.1.4}
F_{\beta_M}^{j_M}F_{\beta_{M-1}}^{j_{M-1}}\ldots
F_{\beta_{M'+1}}^{j_{M'+1}}v_i,\qquad j_k\in\mathbb{Z}_+,\quad
i\in\{1,2,\ldots,d\},
\end{equation}
form a basis of the graded vector space $N(\mathfrak{q}^+,\lambda)$.
\end{proposition}

The homogeneity of the elements \eqref{2.1.4} follows from the fact that
$F_\beta$ are the weight vectors of the $U_q\mathfrak{g}$-module
$U_q\mathfrak{g}$, whose weights are the same as in the classical
case.

Let $U_q\mathfrak{k}\subset U_q\mathfrak{g}$ be the Hopf subalgebra
generated by $\{E_i,F_i\}_{i \ne l_0}\,\cup\, \{K_j^{\pm
1}\}_{j=1,2,\ldots,l}$.

\begin{corollary}\label{1.2.3}
The homogeneous components of the graded vector space
$N(\mathfrak{q}^+,\lambda)$ are finite dimensional weight
$U_q\mathfrak{k}$-modules. The multiplicities of the weights are the same
as in the classical case $q=1$.
\end{corollary}

The algebra $\mathbb{C}[\mathfrak{p}^-]_q$ is a domain, and the homogeneous
component $\mathbb{C}[\mathfrak{p}^-]_{q,1}$ generates
$\mathbb{C}[\mathfrak{p}^-]_q$, see \cite{HeckenbergerKolb04} . Our
immediate intension is to find some defining relations.

The explicit form of the multiplication
$$
m:\mathbb{C}[\mathfrak{p}^-]_q\otimes\mathbb{C}[\mathfrak{p}^-]_q\to
\mathbb{C}[\mathfrak{p}^-]_q,\qquad m:f_1\otimes f_2\mapsto f_1f_2
$$
implies the relation
\begin{equation}\label{3.1.57}
\varphi\psi=
m\check{R}_{\mathbb{C}[\mathfrak{p}^-]_q,\mathbb{C}[\mathfrak{p}^-]_q}
(\varphi\otimes\psi),\qquad\varphi,\psi\in\mathbb{C}[\mathfrak{p}^-]_q.
\end{equation}

In fact, $N(\mathfrak{q}^+,0)$ is an object in the tensor category of
$U_q\mathfrak{g}^\mathrm{cop}$-modules, and the linear maps
$$
\left(\check{R}_{\mathbb{C}[\mathfrak{p}^-]_q,\mathbb{C}[\mathfrak{p}^-]_q}
\right)^*:N(\mathfrak{q}^+,0)\otimes N(\mathfrak{q}^+,0)\to
N(\mathfrak{q}^+,0)\otimes N(\mathfrak{q}^+,0),
$$
$$
\triangle_0:N(\mathfrak{q}^+,0)\to N(\mathfrak{q}^+,0)\otimes
N(\mathfrak{q}^+,0)
$$
are morphisms in this category. On the other hand,
$$
\left(\check{R}_{\mathbb{C}[\mathfrak{p}^-]_q,\mathbb{C}[\mathfrak{p}^-]_q}
\right)^*v(\mathfrak{q}^+,0)\otimes
v(\mathfrak{q}^+,0)=v(\mathfrak{q}^+,0)\otimes v(\mathfrak{q}^+,0),
$$
$$
\triangle_0:v(\mathfrak{q}^+,0)=v(\mathfrak{q}^+,0)\otimes
v(\mathfrak{q}^+,0),
$$
with $v(\mathfrak{q}^+,0)$ being a generator of the
$U_q\mathfrak{g}^\mathrm{cop}$-module $N(\mathfrak{q}^+,0)$. Hence
$\left(\check{R}_{\mathbb{C}[\mathfrak{p}^-]_q,\mathbb{C}[\mathfrak{p}^-]_q}
\right)^*\triangle_0=\triangle_0$. It remains to pass to the adjoint linear
maps. Then the relation \eqref{3.1.57} is to be treated as a commutativity
for $\mathbb{C}[\mathfrak{p}^-]_q$ viewed as an algebra in the braided
tensor category of weight $U_q\mathfrak{b}^-$-finite
$U_q\mathfrak{g}$-modules \cite{Pareigis} .

Denote by $\mathcal{L}$ the kernel of the linear map
$$
\mathbb{C}[\mathfrak{p}^-]_{q,1}\otimes\mathbb{C}[\mathfrak{p}^-]_{q,1}\to
\mathbb{C}[\mathfrak{p}^-]_{q,2},\qquad\varphi\otimes\psi\mapsto\varphi\psi-
m\check{R}_{\mathbb{C}[\mathfrak{p}^-]_q,\mathbb{C}[\mathfrak{p}^-]_q}
(\varphi\otimes\psi).
$$
In the classical limit $q=1$, $\mathcal{L}$ appears to be a subspace of the
antisymmetric tensors.

Let $\mathfrak{k}_\mathrm{ss}=[\mathfrak{k},\mathfrak{k}]$ and
$U_q\mathfrak{k}_\mathrm{ss}\subset U_q\mathfrak{g}$ be the Hopf subalgebra
generated by $\{K_j^{\pm 1}, E_j, F_j\}_{j\ne l_0}$. It follows from
\eqref{3.1.57} that $m\mathcal{L}=0$. To be rephrased, the elements of
$\mathcal{L}$ constitute quadratic relations. We give a description of this
subspace in terms of the morphism of $U_q\mathfrak{k}$-modules
$\widetilde{R}_{\mathbb{C}[\mathfrak{p}^-]_{q,1},
\mathbb{C}[\mathfrak{p}^-]_{q,1}}:\mathbb{C}[\mathfrak{p}^-]_{q,1}\otimes
\mathbb{C}[\mathfrak{p}^-]_{q,1}\to\mathbb{C}[\mathfrak{p}^-]_{q,1}\otimes
\mathbb{C}[\mathfrak{p}^-]_{q,1}$ determined by the relations similar to
\eqref{R-check}, \eqref{A.2.21}, with $\mathfrak{g}$ being replaced by
$\mathfrak{k}_\mathrm{ss}$.

Here we consider various relations between classical and quantum cases. It
is convenient to replace temporarily the ground field
$\mathbb{C}(q^\frac1s)$ by the field of complex numbers, assuming instead
that $q\in(0,1)$. Such $q$ are not the roots of unity.

\begin{proposition}\label{3.1.45} There exists a unique negative eigenvalue
of the linear map $\widetilde{R}_{\mathbb{C}[\mathfrak{p}^-]_{q,1},
\mathbb{C}[\mathfrak{p}^-]_{q,1}}$. This eigenvalue is
$-q^\frac4{(H_0,H_0)}$ and its multiplicity is
$$\dfrac{\dim\mathfrak{p}^-(\dim\mathfrak{p}^--1)}2.$$
\end{proposition}

{\it Proof.} As $\mathbb{C}[\mathfrak{p}^-]_{q,1}$ generates the algebra
$\mathbb{C}[\mathfrak{p}^-]_q$ and $\dim\mathbb{C}[\mathfrak{p}^-]_{q,2}=
\frac{\dim\mathfrak{p}^-(\dim\mathfrak{p}^-+1)}2$, one deduces that
$\dim\mathcal{L}\le\frac{\dim\mathfrak{p}^-(\dim\mathfrak{p}^--1)}2$. Hence
the desired statement is a consequence of the following lemmas.

\begin{lemma}\label{quadratic1}
$\mathcal{L}$ contains all the eigenvectors of
$\widetilde{R}_{\mathbb{C}[\mathfrak{p}^-]_{q,1},
\mathbb{C}[\mathfrak{p}^-]_{q,1}}$ with negative eigenvalues.
\end{lemma}

\begin{lemma}\label{quadratic2}
The dimension of the eigenspace of
$\widetilde{R}_{\mathbb{C}[\mathfrak{p}^-]_{q,1},
\mathbb{C}[\mathfrak{p}^-]_{q,1}}$ with the eigenvalue
$-q^\frac4{(H_0,H_0)}$ is at least
$\frac{\dim\mathfrak{p}^-(\dim\mathfrak{p}^--1)}2$.
\end{lemma}

{\it Proof of Lemma \ref{quadratic1}.} Let $\mathcal{L}'$ be the spectral
subspace of the linear map
$\widetilde{R}_{\mathbb{C}[\mathfrak{p}^-]_{q,1},
\mathbb{C}[\mathfrak{p}^-]_{q,1}}$ associated to the negative half-axis
$(-\infty,0)$. Obviously, $\mathcal{L}$ and $\mathcal{L}'$ are
$U_q\mathfrak{k}$-submodules, so it suffices to prove that
$\mathcal{L}\supset\mathcal{L}'$. In the classical case $q=1$ the
multiplicities of simple weight $U\mathfrak{k}$-modules in
$(\mathfrak{p}^-)^*\otimes(\mathfrak{p}^-)^*$ do not exceed 1 since the
weight subspaces of the $U\mathfrak{k}$-module $\mathfrak{p}^-$ are one
dimensional \cite{Zhelob_engl} . Thus by Propositions \ref{A.2.30},
\ref{100107}, these multiplicities are just 1 in the quantum case as well.
Hence the subspaces $\mathcal{L}$ and $\mathcal{L}'$ are determined by the
respective $U_q\mathfrak{k}$-spectra, i.e., by the sets of highest weights
of their simple $U_q\mathfrak{k}$-submodules.

In the case $q=1$, the sets of highest weights coincide. What remains to do
now is to trace their dependence on $q\in(0,1]$. It follows from
\eqref{A.2.43} that the spectrum of
$\widetilde{R}_{\mathbb{C}[\mathfrak{p}^-]_{q,1},
\mathbb{C}[\mathfrak{p}^-]_{q,1}}$ is on the real axis and does not contain
$0$. Thus the analytic dependence of the spectral projection associated
with the negative half-axis follows from the analytic dependence of
$\widetilde{R}_{\mathbb{C}[\mathfrak{p}^-]_{q,1},
\mathbb{C}[\mathfrak{p}^-]_{q,1}}$ itself.

Now trace the dependence on $q$ of the operators
$$E_j,F_j,\quad j\ne l_0,\qquad H_i,\quad i=1,2,\ldots,l,$$
acting in $\mathcal{L}$, and of the operator
$\widetilde{R}_{\mathbb{C}[\mathfrak{p}^-]_{q,1},
\mathbb{C}[\mathfrak{p}^-]_{q,1}}$. To do this, choose the basis of the
weight vectors in $\mathbb{C}[\mathfrak{p}^-]_q$ dual to the one described
in Proposition \ref{2.1.1}. Thus we get also bases in
$\mathbb{C}[\mathfrak{p}^-]_{q,1}$,
$\mathbb{C}[\mathfrak{p}^-]_{q,1}\otimes\mathbb{C}[\mathfrak{p}^-]_{q,1}$,
$\mathbb{C}[\mathfrak{p}^-]_{q,2}$. The desired statement
$\mathcal{L}\supset\mathcal{L}'$ follows from the fact that the matrix
elements of the operators $m$ and
$\widetilde{R}_{\mathbb{C}[\mathfrak{p}^-]_{q,1},
\mathbb{C}[\mathfrak{p}^-]_{q,1}}$ with respect to the bases given above
depend analytically on $q\in(0,1]$. \hfill $\blacksquare$

\medskip

{\it Proof} of Lemma \ref{quadratic2}. Consider the subspace
\begin{equation}\label{sub_quadr}
\widetilde{\mathcal{L}}=\left\{a\in
N(\mathfrak{q}^+,-\alpha_{l_0})_{-1}\otimes
N(\mathfrak{q}^+,-\alpha_{l_0})_{-1}\left|\:
\check{R}_{N(\mathfrak{q}^+,-\alpha_{l_0}),
N(\mathfrak{q}^+,-\alpha_{l_0})}a=-a\right.\right\}.
\end{equation}
It suffices to prove the inequality
\begin{equation}\label{ineq_quadr}
\dim\widetilde{\mathcal{L}}\ge
\dfrac{\dim\mathfrak{p}^-(\dim\mathfrak{p}^--1)}2.
\end{equation}
In fact, the restrictions of the linear maps
$\check{R}_{N(\mathfrak{q}^+,-\alpha_{l_0}),
N(\mathfrak{q}^+,-\alpha_{l_0})}$,
$\widetilde{R}_{N(\mathfrak{q}^+,-\alpha_{l_0}),
N(\mathfrak{q}^+,-\alpha_{l_0})}$ to the subspace
$N(\mathfrak{q}^+,-\alpha_{l_0})_{-1}\otimes
N(\mathfrak{q}^+,-\alpha_{l_0})_{-1}$ differ only by the scalar multiplier
$q^{-\frac4{(H_0,H_0)}}$ (while comparing the operators of permutation of
the multipliers in the tensor categories of
$U_q\mathfrak{g}^\mathrm{cop}$-modules and
$U_q\mathfrak{k}_\mathrm{ss}^\mathrm{cop}$-modules, we use the reduced
expression for $w_0\in W$ as above).

In order to prove \eqref{ineq_quadr}, consider morphisms of
$U_q\mathfrak{g}^\mathrm{cop}$-modules
$$
N(\mathfrak{q}^+,w\rho-\rho)\to N(\mathfrak{q}^+,-\alpha_{l_0})\otimes
N(\mathfrak{q}^+,-\alpha_{l_0}),\qquad w\in\,
 W^\mathbb{S}\;\&\;l(w)=2,
$$
and the images of homogeneous components
$N(\mathfrak{q}^+,w^{-1}\rho-\rho)_{-2}$. It suffices to prove that the sum
of images has dimension $\dim\mathfrak{p}^-(\dim\mathfrak{p}^--1)/2$, and
the linear map $\check{R}_{N(\mathfrak{q}^+,-\alpha_{l_0}),
N(\mathfrak{q}^+,-\alpha_{l_0})}$ is $-1$ when restricted to each of the
images.

We are going to elaborate the following lemma proved by Kostant
\cite{Kostant_1961} , p. 359 -- 360;\ \ \cite{Rocha} , p. 348, with $r=2$.

\begin{lemma}\label{Kostant-lemma}
Consider the $U\mathfrak{k}$-module $(\mathfrak{p}^-)^{\wedge r}$,
$r=1,2,\ldots,\dim\mathfrak{p}^-$. Its isotypic components are simple
$U\mathfrak{k}$-modules whose weight subspaces are one dimensional, and the
set of weights is
$$
\left\{w\rho-\rho\left|\:w\in\,
 W^\mathbb{S}\;\&\;l(w)=r\right.\right\}.
$$
\end{lemma}

A similar result is also valid for $q\in(0,1)$, as the multiplicities in
decompositions of tensor products remain intact when passing from the
classical case $q=1$ to the quantum case \cite{Jant} . This implies the
desired estimate for the dimension of the sum of homogeneous components
$N(\mathfrak{q}^+,w\rho-\rho)_{-2}$. It remains to show that the linear map
$\check{R}_{N(\mathfrak{q}^+,-\alpha_{l_0}),
N(\mathfrak{q}^+,-\alpha_{l_0})}$ when restricted to an image of any
morphism of $U_q\mathfrak{g}^\mathrm{cop}$-modules
$$
N(\mathfrak{q}^+,w\rho-\rho)\to N(\mathfrak{q}^+,-\alpha_{l_0})\otimes
N(\mathfrak{q}^+,-\alpha_{l_0}),\qquad w\in\,W ^\mathbb{S}\;\&\;l(w)=2,
$$
is $-1$. It suffices to prove that the restriction is $\pm 1$ since it
follows from the continuity in $q$ that we are still inside the spectral
subspace associated with the non-positive part of the spectrum.

Let us take a closer look at \eqref{A.2.43}. It follows from the proof of
the relation expounded in \ \cite{Drinf2_engl}~, p. 239, that the linear
map $\check{R}_{N(\mathfrak{q}^+,\mu),N(\mathfrak{q}^+,\lambda)}
\check{R}_{N(\mathfrak{q}^+,\lambda),N(\mathfrak{q}^+,\mu)}$ when
restricted to the image of a morphism $N(\mathfrak{q}^+,\nu)\to
N(\mathfrak{q}^+,\lambda)\otimes N(\mathfrak{q}^+,\mu)$ is just the scalar
multiplier
$$
q^{-(\mu,\mu+2\rho)-(\lambda,\lambda+2\rho)+(\nu,\nu+2\rho)}=
q^{-(\mu+\rho,\mu+\rho)-(\lambda+\rho,\lambda+\rho)+(\nu+\rho,\nu+\rho)+
(\rho,\rho)}.
$$
Substitute $\lambda=\mu=-\alpha_{l_0}$, $\nu=w\rho-\rho$ to the right-hand
side to get $1$, as the weights $-\alpha_{l_0}-\rho$, $\nu+\rho$, $\rho$
are in the same $W$-orbit and hence have the same length.

\medskip

The proof of Proposition \ref{3.1.45} is finished. \hfill $\blacksquare$

\begin{proposition}\label{3.1.47}
$\mathbb{C}[\mathfrak{p}^-]_q$ is a quadratic algebra with the space of
generators $\mathbb{C}[\mathfrak{p}^-]_{q,1}$ and the space of relations
$\mathcal{L}$.
\end{proposition}

{\it Proof.} Consider the quadratic algebra
$F=T(\mathbb{C}[\mathfrak{p}^-]_{q,1})/(\mathcal{L})$ whose space of
generators is $\mathbb{C}[\mathfrak{p}^-]_{q,1}$ and the space of relations
is $\mathcal{L}$. The natural homomorphism of graded algebras\ \
$\mathscr{I}:\;F\,\to\,\mathbb{C}[\mathfrak{p}^-]_q$ is surjective since
$\mathbb{C}[\mathfrak{p}^-]_{q,1}$ generates
$\mathbb{C}[\mathfrak{p}^-]_q$. The injectivity of $\mathscr{I}$ follows
from the fact that the dimensions of the graded components are the same:
\begin{equation}\label{eq_dim}
\dim\mathbb{C}[\mathfrak{p}^-]_{q,j}=
\binom{j+\dim\mathfrak{p}^--1}{\dim\mathfrak{p}^--1},\qquad\dim F_j=
\binom{j+\dim\mathfrak{p}^--1}{\dim\mathfrak{p}^--1}.
\end{equation}
The first equality in \eqref{eq_dim} can be obtained using the monomial
basis \eqref{2.1.4} of $N(\mathfrak{q}^+,0)$, and the second one is due to
the monomial basis of $F$ described as follows.

Just as in the classical case $q=1$, the weights of the
$U_q\mathfrak{k}$-module $\mathbb{C}[\mathfrak{p}^-]_{q,1}$ are of the form
$-\alpha_{l_0}-\sum\limits_{i\ne l_0}n_i\alpha_i$, and every weight
subspace is one dimensional. Introduce a linear order on the set of weights
of this $U_q\mathfrak{k}$-module corresponding to the lexicographical order
on the set of strings
$(-n_1,-n_2,\cdots,-n_{l_0-1},-1,-n_{l_0+1},\cdots,-n_l)$ formed by the
decomposition coefficients in simple roots. Choose a basis
$\{z_1,z_2,\cdots,z_{\dim\mathfrak{p}^-}\}$ of the weight vectors in
$\mathbb{C}[\mathfrak{p}^-]_{q,1}$, dual to the monomial basis in
$N(\mathfrak{q}^+,0)_{-1}$, and impose an order on its elements
corresponding by the growth of the weights. In view of \eqref{A.2.21}, it
is easy to prove that the tensors
$$
z_i\otimes z_j+q^{-\frac4{(H_\mathbb{S},H_\mathbb{S})}}
\sum_{k<m}\widetilde{R}_{ij}^{km}z_k\otimes z_m,\qquad i>j,
$$
form a non-commutative Gr\"obner basis, hence
\begin{equation}\label{3.1.58}
\left\{\left.z_1^{j_1}z_2^{j_2}z_3^{j_3}\ldots
z_{\dim\mathfrak{p}^-}^{j_{\dim\mathfrak{p}^-}}\right|\:
j_1<j_2<\ldots<j_{\dim\mathfrak{p}^-}\right\}
\end{equation}
form a basis of $F$ \cite{Bergman} .\hfill $\blacksquare$

\begin{remark}\label{Q_q}
In the basis \eqref{3.1.58}, the action of the generators $E_i$, $F_i$,
$K_i^{\pm 1}$ is given by matrices whose elements belong to the field of
rational functions $\mathbb{Q}(q)$ and do not have poles at $q\in(0,1]$, as
it follows from the definitions and Proposition \ref{2.2.2}.
\end{remark}

\section{\sffamily\bfseries\!\!\!\!\!\!. A FIRST ORDER DIFFERENTIAL
CALCULUS}\label{new_FODC}

\indent

A $\mathbb{C}[\mathfrak{p}^-]_q$-bimodule $\Lambda^1(\mathfrak{p}^-)_q$ of
1-forms on the quantum vector space $\mathfrak{p}^-$ was introduced in Sec.
\ref{calculi}, together with the differential
$d:\mathbb{C}[\mathfrak{p}^-]_q\to\Lambda^1(\mathfrak{p}^-)_q$. This
Section presents a description of $\Lambda^1(\mathfrak{p}^-)_q$ in terms of
generators and relations which implies, in particular, that
$\Lambda^1(\mathfrak{p}^-)_q$ is a linear span of the set
$\{f_1df_2f_3|\:f_1,f_2,f_3\in\mathbb{C}[\mathfrak{p}^-]_q\}$.

Let us view the problem in a more general sense. We introduce quantum
analogues of fiberwise linear functions for the homogeneous holomorphic
vector bundles. This way one gets 1-forms in the case of a tangent bundle
for $\mathfrak{p}^-$.

Let $\lambda\in\mathscr{P}_+$. Consider the following morphisms
\begin{align*}
\triangle_{\mathrm{left},\lambda}^+: & \,N(\mathfrak{q}^+,\lambda)\to
N(\mathfrak{q}^+,0)\otimes N(\mathfrak{q}^+,\lambda),
\\ \triangle_{\mathrm{right},\lambda}^+: & \,N(\mathfrak{q}^+,\lambda)\to
N(\mathfrak{q}^+,\lambda)\otimes N(\mathfrak{q}^+,0)
\end{align*}
in the category of $U_q\mathfrak{g}^\mathrm{cop}$-modules, defined by their
action on the generators:
\begin{align*}
\triangle_{\mathrm{left},\lambda}^+: & \,v(\mathfrak{q}^+,\lambda)\to
v(\mathfrak{q}^+,0)\otimes v(\mathfrak{q}^+,\lambda),
\\ \triangle_{\mathrm{right},\lambda}^+: & \,v(\mathfrak{q}^+,\lambda)\to
v(\mathfrak{q}^+,\lambda)\otimes v(\mathfrak{q}^+,0).
\end{align*}

The following relations are immediate consequences of the definitions
\begin{gather}
(\operatorname{id}\otimes\triangle_{\mathrm{left},\lambda}^+)
\triangle_{\mathrm{left},\lambda}^+=(\triangle^+\otimes\operatorname{id})
\triangle_{\mathrm{left},\lambda}^+,\label{4.1.3}
\\ (\triangle_{\mathrm{right},\lambda}^+\otimes\operatorname{id})
\triangle_{\mathrm{right},\lambda}^+=(\operatorname{id}\otimes\triangle^+)
\triangle_{\mathrm{right},\lambda}^+\label{4.1.4}
\\
(\varepsilon^+\otimes\operatorname{id})\triangle_{\mathrm{left},\lambda}^+=
(\operatorname{id}\otimes\varepsilon^+)\triangle_{\mathrm{right},\lambda}^+=
\operatorname{id},\label{4.1.5}
\\ (\operatorname{id}\otimes\triangle_{\mathrm{right},\lambda}^+)
\triangle_{\mathrm{left},\lambda}^+=(\triangle_{\mathrm{left},\lambda}^+
\otimes\operatorname{id})\triangle_{\mathrm{right},\lambda}^+.\label{4.1.6}
\end{gather}
In particular, the latest relation follows from
\begin{eqnarray*}
(\operatorname{id}\otimes\triangle_{\mathrm{right},\lambda}^+)
\triangle_{\mathrm{left},\lambda}^+v(\mathfrak{q}^+,\lambda) &=&
v(\mathfrak{q}^+,0)\otimes v(\mathfrak{q}^+,\lambda)\otimes
v(\mathfrak{q}^+,0),
\\ (\triangle_{\mathrm{left},\lambda}^+\otimes\operatorname{id})
\triangle_{\mathrm{right},\lambda}^+v(\mathfrak{q}^+,\lambda) &=&
v(\mathfrak{q}^+,0)\otimes v(\mathfrak{q}^+,\lambda)\otimes
v(\mathfrak{q}^+,0).
\end{eqnarray*}

Consider the category of graded vector spaces and the space
$\Gamma(\mathfrak{p}^-,\lambda)_q$ which is dual to
$N(\mathfrak{q}^+,\lambda)$ in this category. Equip
$\Gamma(\mathfrak{p}^-,\lambda)_q$ with a structure of
$\mathbb{C}[\mathfrak{p}^-]_q$-bimodule by the duality:
$$
m_{\mathrm{left},\lambda}\overset{\mathrm{def}}{=}
(\triangle_{\mathrm{left},\lambda}^+)^*,\qquad
m_{\mathrm{right},\lambda}\overset{\mathrm{def}}{=}
(\triangle_{\mathrm{right},\lambda}^+)^*.
$$
In a similar way, equip $\Gamma(\mathfrak{p}^-,\lambda)_q$ with a structure
of $U_q\mathfrak{g}$-module using the antipode $S^{-1}$ of the Hopf algebra
$U_q\mathfrak{g}^\mathrm{cop}$.

The relations \eqref{4.1.3} -- \eqref{4.1.6} imply the following statement.

\begin{proposition}
For any $\lambda\in \mathscr{P}_+$ the graded vector space
$\Gamma(\mathfrak{p}^-,\lambda)_q$ is a $U_q\mathfrak{g}$-module
$\mathbb{C}[\mathfrak{p}^-]_q$-bimodule.
\end{proposition}

Let $N(\mathfrak{q}^+,\lambda)_\mathrm{highest}$ be the highest homogeneous
component of the graded vector space $N(\mathfrak{q}^+,\lambda)$ and
$P_\mathrm{highest}$ be the projection in $N(\mathfrak{q}^+,\lambda)$ onto
the subspace $N(\mathfrak{q}^+,\lambda)_\mathrm{highest}$, parallel to the
sum of all other homogeneous components.

\begin{lemma}\label{lemma4}
The linear maps
$$
(\operatorname{id}\otimes
P_\mathrm{highest})\triangle_{\mathrm{left},\lambda}^+:
\,N(\mathfrak{q}^+,\lambda)\to N(\mathfrak{q}^+,0)\otimes
N(\mathfrak{q}^+,\lambda)_\mathrm{highest},
$$
$$
(P_\mathrm{highest}\otimes\operatorname{id})
\triangle_{\mathrm{right},\lambda}^+:\,N(\mathfrak{q}^+,\lambda)\to
N(\mathfrak{q}^+,\lambda)_\mathrm{highest}\otimes N(\mathfrak{q}^+,0)
$$
are injective.
\end{lemma}

{\it Proof.} We restrict ourselves to proving the injectivity of the first
linear map. As for the second one, the proof is similar.

Elaborate the same notation as in the statement of Proposition \ref{2.1.1}
and use the standard order on the set of weights of the
$U_q\mathfrak{k}$-module $N(\mathfrak{q}^+,\lambda)_\mathrm{highest}\cong
L(\mathfrak{k},\lambda)$. The injectivity of the first linear map in
question follows from that proposition and the fact that for any weight
vector $v\in N(\mathfrak{q}^+,\lambda)_\mathrm{highest}$, one has
$$
(\operatorname{id}\otimes
P_\mathrm{highest})\triangle_{\mathrm{left},\lambda}^+
\left(F_{\beta_M}^{j_M}F_{\beta_{M-1}}^{j_{M-1}}\ldots
F_{\beta_{M'+1}}^{j_{M'+1}}v\right)=\mathrm{const}\cdot
F_{\beta_M}^{j_M}F_{\beta_{M-1}}^{j_{M-1}}\ldots
F_{\beta_{M'+1}}^{j_{M'+1}}v(\mathfrak{q}^+,0)\otimes v,
$$
up to the terms possessing lower weights than $v$ in the second tensor
factor. Here $\mathrm{const}\ne 0$ depends on the weight of $v$. \hfill
$\blacksquare$

\begin{proposition}\label{prop14}
\begin{enumerate}
\item The bimodule $\Gamma(\mathfrak{p}^-,\lambda)_q$ over
$\mathbb{C}[\mathfrak{p}^-]_q$ is a free left and a free right
$\mathbb{C}[\mathfrak{p}^-]_q$-module.

\item With $\Gamma(\mathfrak{p}^-,\lambda)_{q,\mathrm{lowest}}$ being the
lowest homogeneous component of $\Gamma(\mathfrak{p}^-,\lambda)_q$, one has
the following isomorphisms of $U_q\mathfrak{k}$-modules:
\begin{align}
\mathbb{C}[\mathfrak{p}^-]_q\otimes
\Gamma(\mathfrak{p}^-,\lambda)_{q,\mathrm{lowest}} & \overset{\approx}{\to}
\Gamma(\mathfrak{p}^-,\lambda)_q, & f\otimes v & \mapsto fv, & \label{4.1.7}
\\ \Gamma(\mathfrak{p}^-,\lambda)_{q,\mathrm{lowest}}\otimes
\mathbb{C}[\mathfrak{p}^-]_q & \overset{\approx}{\to}
\Gamma(\mathfrak{p}^-,\lambda)_q, & v\otimes f & \mapsto vf. & \label{4.1.8}
\end{align}
\end{enumerate}
\end{proposition}

{\it Proof.} The first claim follows from the second one. The morphisms of
$U_q\mathfrak{k}$-modules \eqref{4.1.7}, \eqref{4.1.8} are the morphisms of
graded vector spaces. Proposition \ref{2.1.1} implies by the duality that
the dimensions of the related homogeneous components are finite and equal.
Thus, \eqref{4.1.7}, \eqref{4.1.8} are one-to-one because they are onto,
what is in turn due to Lemma \ref{lemma4}. \hfill $\blacksquare$

\medskip

Of course, $\Gamma(\mathfrak{p}^-,\lambda)_q$ is not a free
$\mathbb{C}[\mathfrak{p}^-]_q$-bimodule. We are to find relations between
the generators from $\Gamma(\mathfrak{p}^-,\lambda)_{q,\mathrm{lowest}}$.
Since $\Gamma(\mathfrak{p}^-,\lambda)_q$ is a lowest weight
$U_q\mathfrak{g}$-module, one has a well defined morphism of
$U_q\mathfrak{g}$-modules
$$
\check{R}_{0,\lambda}\overset{\mathrm{def}}{=}
\check{R}_{\mathbb{C}[\mathfrak{p}^-]_q,\Gamma(\mathfrak{p}^-,\lambda)_q}:
\mathbb{C}[\mathfrak{p}^-]_q\otimes\Gamma(\mathfrak{p}^-,\lambda)_q\to
\Gamma(\mathfrak{p}^-,\lambda)_q\otimes\mathbb{C}[\mathfrak{p}^-]_q.
$$

\begin{proposition}\label{prop15}
\begin{enumerate}
\item $m_{\mathrm{left},\lambda}=
m_{\mathrm{right},\lambda}\cdot\check{R}_{0,\lambda}$,
$m_{\mathrm{right},\mu}=m_{\mathrm{left},\mu}\cdot\check{R}_{\mu,0}$.

\item $\check{R}_{0,\lambda}:\mathbb{C}[\mathfrak{p}^-]_{q,1}\otimes
\Gamma(\mathfrak{p}^-,\lambda)_{q,\mathrm{lowest}}\to
\Gamma(\mathfrak{p}^-,\lambda)_{q,\mathrm{lowest}}\otimes
\mathbb{C}[\mathfrak{p}^-]_{q,1}$.
\end{enumerate}
\end{proposition}

{\it Proof.} The first claim in 1) follows by passing to adjoint operators.
In fact, the claim is equivalent to the coincidence of the morphisms of
$U_q\mathfrak{g}^\mathrm{cop}$-modules
\begin{equation}\label{4.1.9}
\left(\check{R}_{0,\lambda}\right)^*\triangle_{\mathrm{right},\lambda}^+=
\triangle_{\mathrm{left},\lambda}^+.
\end{equation}
In turn, \eqref{4.1.9} follows from the fact that the vector
$v(\mathfrak{q}^+,\lambda)$ generates the
$U_q\mathfrak{g}^\mathrm{cop}$-module $N(\mathfrak{q}^+,\lambda)$, together
taking into account the following relations:
\begin{gather*}
\triangle_{\mathrm{right},\lambda}^+v(\mathfrak{q}^+,\lambda)=
v(\mathfrak{q}^+,\lambda)\otimes v(\mathfrak{q}^+,0),\qquad
\triangle_{\mathrm{left},\lambda}^+v(\mathfrak{q}^+,\lambda)=
v(\mathfrak{q}^+,0)\otimes v(\mathfrak{q}^+,\lambda),
\\ \check{R}_{0,\lambda}^*(v(\mathfrak{q}^+,\lambda)\otimes
v(\mathfrak{q}^+,0))=v(\mathfrak{q}^+,0)\otimes v(\mathfrak{q}^+,\lambda).
\end{gather*}

The second claim in 1) can be proved in a similar way. The second statement
of the Proposition is due to the fact that $\check{R}_{0,\lambda}$ is a
morphism of $U_q\mathfrak{g}$-modules, so it commutes with the linear map
$H_0\otimes 1+1\otimes H_0$. \hfill $\blacksquare$

\begin{corollary}
Let $\{z_i\}$ be a basis of the finite dimensional vector space
$\mathbb{C}[\mathfrak{p}^-]_{q,1}$ and $\{\gamma_i\}$ a basis of the finite
dimensional vector space
$\Gamma(\mathfrak{p}^-,\lambda)_{q,\mathrm{lowest}}$. There exists a unique
 matrix $(\check{R}_{ij}^{km}(\lambda))$ such that
\begin{equation}\label{4.1.10}
\check{R}_{0,\lambda}(z_i\otimes\gamma_j)=
\sum_{k,m}\check{R}_{ij}^{km}(\lambda)\gamma_k\otimes z_m,
\end{equation}
with\ \ $i,m\in \{1,2,\ldots,\dim\mathbb{C}[\mathfrak{p}^-]_{q,1}\}$,\ \
$j,k\in
\{1,2,\ldots,\dim\Gamma(\mathfrak{p}^-,\lambda)_{q,\mathrm{lowest}}\}$.
\end{corollary}

\begin{remark}
To find the functions $\check{R}_{ij}^{km}(\lambda)$ seems to be an
intricate problem, as their definition involves actions of the universal
R-matrix in the tensor products of infinite dimensional
$U_q\mathfrak{g}$-modules. We show that actually this is not a problem.
Just as in Sec. \ref{new_item}, consider the Hopf subalgebra
$U_q\mathfrak{k}_\mathrm{ss}\subset U_q\mathfrak{g}$ generated by $K_i^{\pm
1}$, $E_i$, $F_i$ with $i\ne l_0$. The action of the universal R-matrix of
the Hopf algebra $U_q\mathfrak{g}$ on the vectors from
$\mathbb{C}[\mathfrak{p}^-]_{q,1}\otimes
\Gamma(\mathfrak{p}^-,\lambda)_{q,\mathrm{lowest}}$ differs from the action
of the universal R-matrix of the Hopf algebra $U_q\mathfrak{k}_\mathrm{ss}$
only by a multiplier
\begin{equation}\label{4.1.11}
\mathrm{const}=q_{l_0}^\frac{(\lambda,\overline{\omega}_{l_0})}
{(\overline{\omega}_{l_0},\overline{\omega}_{l_0})}.
\end{equation}
To prove \eqref{4.1.11}, it suffices to use \eqref{A.2.21}, the reduced
expression of $w_0$ as in Sec. \ref{new_item}, and the relation
$$
(\alpha_{l_0},\lambda)=(\alpha_{l_0}|_{\mathfrak{h}\cap\mathfrak{k}},
\lambda|_{\mathfrak{h}\cap\mathfrak{k}})+\frac{(\alpha_{l_0},
\overline{\omega}_{l_0})(\overline{\omega}_{l_0},\lambda)}
{(\overline{\omega}_{l_0},\overline{\omega}_{l_0})},\qquad
(\alpha_{l_0},\overline{\omega}_{l_0})=d_{l_0},
$$
which allow to compare the Cartan multipliers $q^{-t_0}$.
\end{remark}

\medskip

The next statement follows from Propositions \ref{prop14}, \ref{prop15}.

\begin{proposition}
The set $\{\gamma_j\}_{j=1,2,\ldots,
\dim\Gamma(\mathfrak{p}^-,\lambda)_{q,\mathrm{lowest}}}$ generates the
$\mathbb{C}[\mathfrak{p}^-]_q$-bimodule $\Gamma(\mathfrak{p}^-,\lambda)_q$,
and
\begin{equation}\label{4.1.13}
z_i\gamma_j=\sum_{k,m}\check{R}_{ij}^{km}(\lambda)\,\gamma_kz_m
\end{equation}
are defining relations.
\end{proposition}

\begin{corollary}
For all $i,j\in \{1,2,\ldots,\dim\mathfrak{p}^-\}$ one has
\begin{equation}\label{4.1.14}
z_idz_j=\sum_{k,m=1}^{\dim\mathfrak{p}^-} \check{R}_{ij}^{km}\,dz_kz_m,
\end{equation}
with $\check{R}_{ij}^{km}=\check{R}_{ij}^{km}(-\alpha_{l_0})$, which
constitutes a defining list of relations between the generators of the
$\mathbb{C}[\mathfrak{p}^-]_q$-bimodule $\Lambda^1(\mathfrak{p}^-)_q$.
\end{corollary}

\section{\sffamily\bfseries\!\!\!\!\!\!. A UNIVERSAL DIFFERENTIAL
CALCULUS}\label{new_diff_univ}

\indent

As it was mentioned in Sec. \ref{calculi}, every first order differential
calculus determines a universal differential calculus over the algebra $F$.
In the present section we give several examples illustrating this
situation.

\begin{example}(No relations.)
Consider a vector space $V$ and a differential calculus over its tensor
algebra $T(V)$. Choose a vector space $V'$ being isomorphic to $V$,
together with an isomorphism $d:V\to V'$. Equip the tensor algebra
$\Omega'=T(V\oplus V')$ with a grading as follows:
$$\deg(v) =0,\quad v\in V;\qquad\deg(v')=1,\quad v'\in V'.$$
Obviously, $\Omega'=\bigoplus\limits_{j\in\mathbb{Z}_+}\Omega_j'$, with
$\Omega_j'=\{t\in\Omega'|\:\deg t=j\}$. Define a linear operator $d'$ in
$\Omega'$ recursively: $d'1=0$,
\begin{gather*}
d'v=dv,\quad v\in V;\qquad d'v'=0,\quad v'\in V',
\\ d'(v\otimes t)=dv\otimes t+v\otimes d't,\qquad v\in V,\;t\in\Omega',
\\ d'(v'\otimes t)=-v'\otimes d't,\qquad v'\in V',\;t\in\Omega'.
\end{gather*}
The pair $(\Omega_1',d'|_{\Omega_0'})$ is a first order differential
calculus over the tensor algebra $T(V)$, and the pair $(\Omega',d')$ is its
universal differential calculus.
\end{example}

\begin{example}(Relations between the coordinates are imposed.)
Consider the two-sided ideal $J_0=\bigoplus\limits_{j\ge 2}(J_0\cap
V^{\otimes j})$ of the tensor algebra
$T(V)=\bigoplus\limits_{j\in\mathbb{Z}_+}V^{\otimes j}$. We show how to
define a differential calculus over the algebra $F=T(V)/J_0$. Let $J_F$ be
the two-sided ideal of $\Omega'$ generated by $J_0$ and $d'J_0$. The
algebra $\Omega^F=\Omega'/J_F$ inherits the grading \
$\Omega^F=\bigoplus\limits_{j\in\mathbb{Z}_+}\Omega_j^F$.

Obviously, $V\hookrightarrow F$. As $d'J_F\subset J_F$, the differential
$d'$ transfers onto $\Omega^F=\Omega'/J_F$. Thus we get a linear map $d_F$
and a differential calculus $(\Omega^F,d^F)$ over the algebra
$F=\Omega_0^F$. It is the universal differential calculus of the first
order differential calculus $(\Omega_1^F,d^F|_F)$, and is well-known in the
quantum group theory \cite{KlSch} , p. 462.
\end{example}

\begin{example}\label{4.1.18} (Relations between the coordinates and the
differentials are imposed.)

We have described above the free first order differential calculus over the
algebra $F=T(V)/J_F$. Now turn to a more realistic example by introducing
$R$-matrix commutation relations between the elements $v\in V$ and $v'\in
V'$. Consider an invertible linear map $\check{\mathcal{R}}:V\otimes V'\to
V'\otimes V$. Let $J_1$ be the subbimodule of the $F$-bimodule $\Omega_1^F$
generated by
\begin{equation}\label{4.1.16}
\{vv'-v'v|\:v\otimes v'=\check{\mathcal{R}}(v\otimes v'),\quad v\in
V,\;v'\in V'\}\subset\Omega_1^F.
\end{equation}
The factorization of $(\Omega_1^F,d^F|_F)$ by $J_1$ leads to a first order
differential calculus over $F$. To obtain the associated universal
differential calculus, consider the homogeneous two-sided ideal $J_M$ of
$\Omega^F$ generated by $J_1$, $d^FJ_1$, and divide out $\Omega^F$ by
$J_M$. It is possible to transfer the differential onto $\Omega^F/J_M$
since $d^FJ_M\subset J_M$.
\end{example}

Sec. \ref{calculi}, \ref{new_FODC} provide a first order differential
calculus $(\Lambda^1(\mathfrak{p}^-)_q,d)$ over the algebra
$\mathbb{C}[\mathfrak{p}^-]_q$, together with the associated universal
differential calculus $(\Lambda(\mathfrak{p}^-)_q,d)$.

It is easy to obtain its description in terms of generators and relations,
in view of the results of Sec. \ref{new_FODC} and Example \ref{4.1.18}. For
that, it suffices to complete the list of relations for the
$\mathbb{C}[\mathfrak{p}^-]_q$-bimodule $\Lambda^1(\mathfrak{p}^-)_q$ with
those deduced from \eqref{4.1.14} via differentiation:
\begin{equation}\label{4.1.17}
dz_i\wedge
dz_j=-\sum_{k,m=1}^{\dim\mathfrak{p}^-}\check{R}_{ij}^{km}\,dz_k\wedge
dz_m.
\end{equation}

It is easy to prove that $\Lambda^j(\mathfrak{p}^-)_q$ is a free left
$\mathbb{C}[\mathfrak{p}^-]_q$-module of rank
$\dbinom{\dim\mathfrak{p}^-}{j}$ and also a free right
$\mathbb{C}[\mathfrak{p}^-]_q$-module of rank
$\dbinom{\dim\mathfrak{p}^-}{j}$, just as in the classical case $q=1$.

Now notice that the linear span of the $j$-forms
$$
dz_{i_1}\wedge dz_{i_2}\wedge\ldots\wedge dz_{i_j},\qquad
i_1,i_2,\ldots,i_j\in\{1,2,\ldots,\dim\mathfrak{p}^-\},
$$
with constant coefficients, is a $U_q\mathfrak{k}$-module. This linear span
will be denoted by $\Lambda^j(\mathfrak{p}^-)_q^\mathrm{const}$.

\begin{proposition}
There exist the isomorphisms of $U_q\mathfrak{k}$-modules
\begin{align}
\mathbb{C}[\mathfrak{p}^-]_q\otimes
\Lambda^j(\mathfrak{p}^-)_q^\mathrm{const} &
\overset{\approx}{\to}\Lambda^j(\mathfrak{p}^-)_q, & f\otimes\omega &
\mapsto f\omega, & \label{4.1.20}
\\
\Lambda^j(\mathfrak{p}^-)_q^\mathrm{const}\otimes
\mathbb{C}[\mathfrak{p}^-]_q & \overset{\approx}{\to}
\Lambda^j(\mathfrak{p}^-)_q, & \omega\otimes f & \mapsto \omega f. &
\label{4.1.21}
\end{align}
\end{proposition}

{\it Proof.} It suffices to use Proposition \ref{prop14} and the definition
of $\Lambda(\mathfrak{p}^-)_q$. \hfill $\blacksquare$

\medskip

We turn to a computation of
$\dim\Lambda^j(\mathfrak{p}^-)_q^\mathrm{const}$.

\begin{lemma}\label{leq}
The dimension of the vector space
$\Lambda^2(\mathfrak{p}^-)_q^\mathrm{const}$ does not exceed its value in
the classical case $q=1$:
\begin{equation}\label{4.1.22}
\dim\Lambda^2(\mathfrak{p}^-)_q^\mathrm{const}\le
\frac{\dim\mathfrak{p}^-(\dim\mathfrak{p}^--1)}2.
\end{equation}
\end{lemma}

{\it Proof.} While proving \eqref{4.1.22}, we may replace the ground field
$\mathbb{C}(q^\frac1s)$ with the ground field $\mathbb{C}$ assuming
$q\in(0,1)$ to be transcendental. Consider the linear map
$$
\widetilde{R}:dz_i\otimes
dz_j\mapsto\sum_{k,m=1}^{\dim\mathfrak{p}^-}\check{R}_{ij}^{km}dz_k\otimes
dz_m
$$
in $\Lambda^1(\mathfrak{p}^-)_q^\mathrm{const}\otimes
\Lambda^1(\mathfrak{p}^-)_q^\mathrm{const}$. One has a natural isomorphism
$$
\Lambda^2(\mathfrak{p}^-)_q^\mathrm{const}\cong
\left\{v\in\Lambda^1(\mathfrak{p}^-)_q^\mathrm{const}
\otimes\Lambda^1(\mathfrak{p}^-)_q^\mathrm{const}\:\left|\:
\widetilde{R}v=-v\right.\right\}.
$$

It is easy to see that all the eigenvalues of the linear map
$\widetilde{R}$ are real and non-zero. A proof reduces to the replacement
of the $\check{R}$-matrix of the Hopf algebra $U_q\mathfrak{g}$ by the
$\widetilde{R}$-matrix of the Hopf algebra $U_q\mathfrak{k}_\mathrm{ss}$,
with the subsequent application of \eqref{A.2.43}. One may assume a basis
$\{z_i\}$ of the vector space $\mathbb{C}[\mathfrak{p}^-]_{q,1}$ to be
chosen as in Sec. \ref{new_item}. Under this choice, the functions
$\check{R}_{ij}^{km}$ in the indeterminate $q$ are analytic on the
half-interval $(0,1]$, hence are continuous for $0<q\le 1$. Thus for all
$q\in(0,1]$ the number of negative eigenvalues of $\widetilde{R}$ according
to their multiplicities is
$\dfrac{\dim\mathfrak{p}^-(\dim\mathfrak{p}^--1)}2$. Hence the dimension of
the eigenspace corresponding to the eigenvalue $-1$ does not exceed
$\dfrac{\dim\mathfrak{p}^-(\dim\mathfrak{p}^--1)}2$. \hfill $\blacksquare$

\medskip

The inequality, that is converse to \eqref{4.1.22}, will be established as
soon as $-1$ is proved to be an eigenvalue of $\widetilde{R}$ with the
multiplicity at least $\dfrac{\dim\mathfrak{p}^-(\dim\mathfrak{p}^--1)}2$.
It suffices to get a similar estimate for the adjoint linear map. The
latter may be identified as the restriction of
$\check{R}_{N(\mathfrak{q}^+,-\alpha_{l_0}),
N(\mathfrak{q}^+,-\alpha_{l_0})}$ to the tensor product of the highest
homogeneous components $N(\mathfrak{q}^+,-\alpha_{l_0})_{-1}\otimes
N(\mathfrak{q}^+,-\alpha_{l_0})_{-1}$. The desired inequality follows from
\eqref{ineq_quadr}. Thus we have proved

\begin{lemma}\label{geq}
The dimension of the vector space
$\Lambda^2(\mathfrak{p}^-)_q^\mathrm{const}$ is at least its value for
$q=1$:
\begin{equation}\label{4.1.34}
\dim\Lambda^2(\mathfrak{p}^-)_q^\mathrm{const}\ge
\frac{\dim\mathfrak{p}^-(\dim\mathfrak{p}^--1)}2.
\end{equation}
\end{lemma}

Now we turn to the higher order differential forms.

\begin{lemma}\label{3.6.27}
For all $j\ge 3$
\begin{equation}\label{4.1.35}
\dim\Lambda^j(\mathfrak{p}^-)_q^\mathrm{const}\le
\binom{\dim\mathfrak{p}^-}{j}.
\end{equation}
\end{lemma}

{\it Proof.} Consider the basis of weight vectors
$\{z_1,z_2,\ldots,z_{\dim\mathfrak{p}^-}\}$ of the $U_q\mathfrak{k}$-module
$\mathbb{C}[\mathfrak{p}^-]_{q,1}$ formed in the final part of Sec.
\ref{new_item}. Impose the lexicographic order relation on the set
$\{dz_i\otimes dz_k\}_{i,k=1,2,\ldots,\dim\mathfrak{p}^-}$. The action of
the universal R-matrix of the Hopf algebra $U_q\mathfrak{k}$ with respect
to this basis is given by a triangular matrix with positive entries on the
principal diagonal (see \eqref{A.2.21}). Thus only the terms $dz_k\otimes
dz_m$ with $k\le m$ contribute to the right-hand side of
\begin{equation}\label{4.1.36}
dz_i\otimes
dz_j=-\sum_{k,m=1}^{\dim\mathfrak{p}^-}\check{R}_{ij}^{km}dz_k\otimes
dz_m,\qquad i\ge j,
\end{equation}
and every element of $\Lambda^j(\mathfrak{p}^-)_q^\mathrm{const}$ belongs
to the linear span of
$$
dz_{i_1}\wedge dz_{i_2}\wedge\cdots\wedge dz_{i_j},\qquad 1\le
i_1<i_2<\ldots<i_j\le\dim\mathfrak{p}^-.\eqno\blacksquare
$$

\medskip

\begin{remark}
The distinguished basis
$$
\{dz_1,dz_2,\ldots,dz_{\dim\mathfrak{p}^-}\}\subset
\Lambda^1(\mathfrak{p}^-)_q^\mathrm{const}
$$
constructed in the proof of Lemma \ref{3.6.27} can be used to identify the
tensor algebra $T((\mathfrak{p}^-)^*)$ with the free non-commutative
algebra $\mathbb{C}\langle
dz_1,dz_2,\ldots,dz_{\dim\mathfrak{p}^-}\rangle$. Consider the two-sided
ideal $I$ of the free algebra generated by the set $G$ of the differences
between the left- and the right-hand sides of \eqref{4.1.36}. The latter
relations work as ``substitution rules'' in the final part of the proof of
Lemma \ref{3.6.27}: the left-hand side, wherever it occurs, is to be
replaced by the right-hand side.
\end{remark}

\medskip

A proof of the following lemma uses the same arguments as in the proof of a
similar result from the work by Heckenberger and Kolb on the de Rham
complexes \cite{HeckenbergerKolb03}~.

\begin{lemma}\label{j_geq_3}
For all $j\ge 3$
\begin{equation}\label{4.1.37}
\dim\Lambda^j(\mathfrak{p}^-)_q^\mathrm{const}\ge
\binom{\dim\mathfrak{p}^-}{j}.
\end{equation}
\end{lemma}

{\it Proof.} If \eqref{4.1.37} holds for $j=3$, it follows from the diamond
lemma \cite{Bergman} that $G$ is a Gr\"obner basis for the two-sided ideal
$I$. In this context, the desired inequality holds also for all $j\ge 3$.
This means that we may restrict ourselves to the special case $j=3$ while
proving Lemma \ref{j_geq_3}.

We identify $(\mathfrak{p}^-)^*$ to
$\Lambda^1(\mathfrak{p}^-)_q^\mathrm{const}$ and introduce the notation
$(\mathfrak{p}^-)^{*\wedge 2}$ for the subspace
$\left\{v\in(\mathfrak{p}^-)^*\otimes(\mathfrak{p}^-)^*\left|\,
\widetilde{R}v=-v\right.\right\}$, cf. \eqref{4.1.36}. It follows from
Lemmas \ref{leq} and \ref{geq} that
\begin{equation}\label{4.1.38}
\dim(\mathfrak{p}^-)^{*\wedge
2}=\frac{\dim\mathfrak{p}^-(\dim\mathfrak{p}^--1)}2.
\end{equation}

Consider the subspaces $L_1=(\mathfrak{p}^-)^{*\wedge
2}\otimes(\mathfrak{p}^-)^*$,
$L_2=(\mathfrak{p}^-)^*\otimes(\mathfrak{p}^-)^{*\wedge 2}$ of the vector
space $(\mathfrak{p}^-)^{*\otimes 3}$, together with a complex of linear
maps
\begin{gather*}
0\to L_1\cap L_2\to L_1\oplus L_2\overset{j}{\to}(\mathfrak{p}^-)^{*\otimes
3}\to\mathbb{C}[\mathfrak{p}^-]_{q,3}\to 0,
\\ j:v_1\oplus v_2\mapsto v_1-v_2,\qquad v_j\in
L_j\subset(\mathfrak{p}^-)^{*\otimes 3},
\end{gather*}
which is exact in all terms except, possibly, $(\mathfrak{p}^-)^{*\otimes
3}$. After computing the Euler characteristic for this complex, we deduce
that
$$
-\dim(L_1\cap L_2)+\dim(L_1\oplus L_2)-\dim\left((\mathfrak{p}^-)^{*\otimes
3}\right)+\dim\mathbb{C}[\mathfrak{p}^-]_{q,3}\le 0.
$$
Now use \eqref{4.1.38}, \eqref{eq_dim} to get
$$
\dim(L_1\cap L_2)\ge n^2(n-1)-n^3+\frac{n(n+1)(n+2)}6=\frac{n(n-1)(n-2)}6,
$$
with $n=\dim\mathfrak{p}^-=\dim(\mathfrak{p}^-)^*$. One deduces from a
description of the universal differential calculus as in Example
\ref{4.1.18} that
$$
\dim\Lambda^3(\mathfrak{p}^-)_q^\mathrm{const}\ge\dim(L_1\cap
L_2)\ge\binom{\dim\mathfrak{p}^-}3.\eqno\blacksquare
$$

\medskip

The next statement follows from the Lemmas proved in this Section.

\begin{proposition}
The homogeneous components $\Lambda^j(\mathfrak{p}^-)_q$ of the
differential graded algebra $\Lambda(\mathfrak{p}^-)_q$ vanish for
$j>\dim\mathfrak{p}^-$. Each of those is a free left and also a free right
$\mathbb{C}[\mathfrak{p}^-]_q$-module of rank
$\dbinom{\dim\mathfrak{p}^-}j$.
\end{proposition}

The following statement is well known in the classical case $q=1$ and can
be proved by reduction to that case.

\begin{proposition}
In $\Lambda(\mathfrak{p}^-)_q$ one has
$\operatorname{Ker}d=\operatorname{Im}d$.
\end{proposition}

{\it Proof.} Consider the introduced above monomial bases \eqref{3.1.58} in
$\mathbb{C}[\mathfrak{p}^-]_q$ and $$\{dz_{i_1}\wedge
dz_{i_2}\wedge\cdots\wedge dz_{i_k}|\:k\in\mathbb{Z}_+,\quad 1\le
i_1<i_2<\ldots<i_k\le\dim\mathfrak{p}^-\}$$ in
$\Lambda(\mathfrak{p}^-)_q^\mathrm{const}$. Use these bases, together with
the isomorphism of vector spaces
$\Lambda(\mathfrak{p}^-)_q\cong\mathbb{C}[\mathfrak{p}^-]_q\otimes
\Lambda(\mathfrak{p}^-)_q^\mathrm{const}$ to get a monomial basis in
$\Lambda(\mathfrak{p}^-)_q$.

In this basis, the matrix elements of $d$ are rational functions from
$\mathbb{Q}(q)$ and do not have poles in $(0,1]$, see \eqref{A.2.21} and
Proposition \ref{2.2.2}. The relation
$\operatorname{Ker}d=\operatorname{Im}d$ at transcendental $q$ follows from
the fact that it holds for $q=1$. Actually, equip the algebra
$\Lambda(\mathfrak{p}^-)_q$ with the grading $\deg(z_j)=\deg(dz_j)=1$,
$j=1,2,\ldots,\dim\mathfrak{p}^-$. Observe that $d$ preserves the
homogeneity degree of differential forms, and the homogeneous components of
$\Lambda(\mathfrak{p}^-)_q$ are finite dimensional. It remains to use the
relation $d^2=0$ and the following well-known result.

Let $A(q)$ be a matrix with the entries from $\mathbb{Q}(q)$ which
satisfies $A(q)^2=0$. Then the associated operator function satisfies
$\dim\operatorname{Ker}A(q)=\dim\operatorname{Im}A(q)$ on a Zariski open
subset. In fact, the function
$\dim\operatorname{Ker}A(q)-\dim\operatorname{Im}A(q)$ takes values in
$\mathbb{Z}_+$ and is upper semicontinuous since both
$\dim\operatorname{Ker}A(q)$ and $-\dim\operatorname{Im}A(q)$ are upper
semicontinuous. \hfill $\blacksquare$

\section{\sffamily\bfseries\!\!\!\!\!\!. A BGG RESOLUTION OF THE TRIVIAL
\boldmath $U_q\mathfrak{g}$-MODULE}\label{BGG-section}

\indent

In this Section we remind the well-known results \ \cite{FeFr,Malikov} .
Let $\lambda=\{\lambda_1,\lambda_2,\ldots,\lambda_l\}\in P$. The Verma
module $M(\lambda)$ with the highest weight $\lambda$ is generated by the
highest weight vector $v(\lambda)$. The defining relations are as follows:
$$
E_iv(\lambda)=0,\qquad K_i^{\pm
1}v(\lambda)=q_i^{\pm\lambda_i}v(\lambda),\quad i=1,2,\ldots,l.
$$

The affine action of the Weyl group $W$ is introduced by
$w\cdot\lambda=w(\lambda+\rho)-\rho$, with $\rho$ being a half-sum of
positive roots.

A non-zero vector $v\in M(\lambda)$ is said to be singular if it is a
weight vector and $E_iv=0$ for all $i=1,2,\ldots,l$. There exists a
one-to-one correspondence between the singular vectors of the weight $\mu$
in $M(\lambda)$ and the non-zero morphisms of Verma modules $M(\mu)\to
M(\lambda)$. Every non-zero morphism of Verma modules is injective and
$\dim\operatorname{Hom}_{U_q\mathfrak{g}}(M(\mu),M(\lambda))\le 1$ since a
similar result is valid in the classical case $q=1$. The details of the
argument can be found in \ \cite{FeFr} , Sec. 4.5.

The next claim is well known in the classical case $q=1$ and can be easily
deduced from the definitions.

\begin{lemma} {\rm(cf. \cite{Dix_engl} , Proposition 7.1.15)}
If $\lambda=(\lambda_1,\lambda_2,\ldots,\lambda_l)\in\mathbb{R}^l$ and
$\lambda_i\in\mathbb{Z}_+$ for some $i\in\{1,2,\ldots,l\}$, then the vector
$v_i=F_i^{\lambda_i+1}v(\lambda)\in M(\lambda)$ is a singular vector of the
weight $s_i\cdot\lambda$.
\end{lemma}

An application of the lemma, together with the argument used
in\cite{Dix_engl}~, yields

\begin{proposition} {\rm(\cite{Dix_engl} , Proposition 7.6.8)}
For any $w\in W$
$$\dim\operatorname{Hom}_{U_q\mathfrak{g}}(M(w\cdot 0),M(0))=1.$$
\end{proposition}

In what follows we fix the embeddings $i_w:M(w\cdot 0)\hookrightarrow
M(0)$, that is, singular vectors which are the images of $v(w\cdot 0)$
under $i_w$. We choose the singular vectors so that the coefficients of
their decomposition with respect to the basis
$$
F_{\beta_1}^{j_1}F_{\beta_2}^{j_2}\ldots F_{\beta_M}^{j_M}v(0),\qquad
j_1,j_2,\ldots j_M\in\mathbb{Z}_+,
$$
belong to the field $\mathbb{Q}(q)$ and do not have poles at $q=1$. The
formulas for the singular vectors of Verma modules were obtained in
\cite{Malikov, Feigin_Malikov, Iohara_Malikov} and in
\cite{Dobrev_singular} ; the formulas for the projections onto subspaces of
singular vectors with fixed weights were obtained by Tolstoy
\cite{Tolstoy}~.

As it is clear from \ \cite{FeFr,Malikov}~, the BGG resolution of the
trivial $U_q\mathfrak{g}$-module $\mathbb{C}$ has the same form as in the
classical case $q=1$:
\begin{equation}\label{8.2.1}
\ldots\overset{d_2}{\longrightarrow}C_1\overset{d_1}{\longrightarrow}C_0
\overset{\epsilon}{\longrightarrow}\mathbb{C}\longrightarrow 0,\qquad
C_j=\bigoplus_{\{w\in W|\:l(w)=j\}}M(w\cdot 0).
\end{equation}
Here $\epsilon:M(0)\to\mathbb{C}$, $\epsilon:v(0)\mapsto 1$ is an obvious
morphism of $U_q\mathfrak{g}$-modules.

The construction of differentials $d_j$ elaborates a partial order on the
Weyl group $W$, the Bruhat order. Recall a definition from \cite{Hum2} .

Consider the oriented graph $\mathbb{G}$ whose vertices are the elements of
$W$ and the edges are such ordered pairs $w'\to w''$ of vertices that
$l(w'')=l(w')+1$, and $w''=w's_\gamma$ with $s_\gamma$ being a reflection
corresponding to a root $\gamma\in\Phi$. By definition, $w'\le w''$ iff
there exists a path from $w'$ to $w''$. This partial order is called the
Bruhat order.

The above order remains intact if in its definition one replaces
$w''=w's_\gamma$ by $w''=s_\gamma w'$ \cite{Hum2} , p. 119. Also, if in the
definition of the graph $\mathbb{G}$ one replaces the equality
$l(w'')=l(w')+1$ by the inequality $l(w'')>l(w')$, one gets a different but
equivalent definition of the Bruhat order \cite{Hum2} , p. 118, 122. We
present a description of the Bruhat order in terms of reduced expressions
for the elements of the Weyl group $W$.

\begin{proposition} (\cite{Hum2} , p. 120)
Let
\begin{equation}\label{8.2.2}
w=s_1s_2\cdots s_{l(w)}
\end{equation}
be a reduced expression of an element $w\in W$. The set $\{w'\in W|\:w'\le
w\;\&\;w'\ne w\}$ coincides with the set of the elements produced by
omitting some (possibly all) multipliers in the right-hand side of
\eqref{8.2.2}:
$$
w'=s_{i_1}s_{i_2}\cdots s_{i_r},\qquad 1\le i(1)<i(2)<\ldots<i(r)\le l(w).
$$
\end{proposition}

The next well-known result provides a correspondence between the Bruhat
order on $W$ and the standard order relation on the subset $\{w\cdot
0|\:w\in W\}$ of the weight lattice $P\cong\mathbb{Z}^l$.

\begin{lemma} [cf. \ \cite{Dix_engl} , proposition 7.7.2]
Let $w\in W$, $\gamma\in\Phi$. Then

i) $(s_\gamma w)\cdot 0=w\cdot 0-n\gamma$, with $n$ a non-zero integer;

ii) if $l(s_\gamma w)>l(w)$, then $(s_\gamma w)\cdot 0<w\cdot 0$;

iii) if $l(s_\gamma w)<l(w)$, then $(s_\gamma w)\cdot 0>w\cdot 0$.
\end{lemma}

Consider the morphisms of $U_q\mathfrak{g}$-modules
$$i_{w_1,w_2}:M(w_1\cdot 0)\to M(w_2\cdot 0),\qquad w_1,w_2\in W,$$
such that $i_{w_2}i_{w_1,w_2}=i_{w_1}$.

The existence and the uniqueness of $i_{w_1,w_2}$ under the assumption
$i_{w_1}M(w_1\cdot 0)\subset i_{w_2}M(w_2\cdot 0)$ are obvious. The
inclusion holds iff $w_1\ge w_2$. This fact can be established in the same
way as in the classical case $q=1$: it suffices to repeat the proofs of \
\cite{Dix_engl} , Lemmas 7.6.10, 7.6.11, with the reference to Lemma 7.6.9
replaced by a reference to the following statement.

\begin{lemma}
\begin{enumerate}
\item For any $i\in\{1,2,\ldots,l\}$, $u\in U_q\mathfrak{n}^-$ there exists
such $N\in\mathbb{N}$ that $F_i^Nu\in U_q\mathfrak{n}^-\cdot F_i$.
\item For any $i\in\{1,2,\ldots,l\}$, $u\in U_q\mathfrak{n}^-$ there exists
such $N\in\mathbb{N}$ that $uF_i^N\in F_i\cdot U_q\mathfrak{n}^-$.
\end{enumerate}
\end{lemma}

{\it Proof.} We start proving the first claim. It suffices to prove it
in the special case $u=F_j$ with $j\in\{1,2,\ldots,l\}$, since $F_j$
generate $U_q\mathfrak{n}^-$. Even more, one may assume $j\ne i$. It
follows from the defining relations between $F_i$, $F_j$ in
$U_q\mathfrak{g}$ that $F_i^{1-a_{ij}}F_j$ is a linear combination of
$F_i^{1-a_{ij}-k}F_jF_i^k$ with $1\le k\le 1-a_{ij}$. Set $N=1-a_{ij}$.

The second claim can be proved in a similar way. \hfill $\blacksquare$

\medskip

This lemma means that the multiplicative subset $F_i^{\mathbb{Z}_+}$ of
$U_q\mathfrak{n}^-$ satisfies both the right and the left Ore conditions.

The following result of \cite{BGG} is crucial.

\begin{lemma}\label{BGG-square}
Consider the oriented graph $\mathbb{G}$ introduced above.
\begin{enumerate}
\item Let $w,w''\in W$, $w\le w''$, and $l(w'')=l(w)+2$. Then the
    number of such $w'\in W$ that there exist the edges $w\to w'$,
    $w'\to w''$, is either zero or two (in the latter case one has a
    quadruple of the elements of the Weyl groups $W$, called a square).

\item It is possible to associate a number $\epsilon(w,w')=\pm 1$ with
    every edge $w'\to w$ in such a way that the product of the numbers
    corresponding to the edges of each square is $-1$.
\end{enumerate}
\end{lemma}

A differential $d_j:C_j\to C_{j-1}$ is defined by
\begin{equation}\label{8.2.3}
d_j|_{M(w\cdot 0)}=\bigoplus_{\{w'\in W|\:w'\to w\}}\epsilon(w,w')i_{w,w'}.
\end{equation}

Obviously, $d_j\circ d_{j-1}=0$ for all $j\in\mathbb{N}$. Moreover, the
submodule $\operatorname{Ker}\epsilon$ of the Verma module $M(0)$ is
generated by $F_iv(0)$, $i=1,2,\ldots,l$. Thus
$\operatorname{Im}d_1=\operatorname{Ker}\epsilon$, and \eqref{8.2.3} is a
complex which is extended by $\epsilon$ in the category of
$U_q\mathfrak{g}$-modules. Its exactness for a transcendental $q$ follows
from the result of Bernstein-Gelfand-Gelfand on the exactness of a similar
complex in the classical case $q=1$. In fact, the weight subspaces of Verma
modules are finite dimensional, and the matrix elements of the linear maps
$d_j|_{M(w\cdot 0)}$ in the bases introduced in Sec. \ref{PBW-R} belong to
the field $\mathbb{Q}(q)$ of rational functions and do not have poles at
$q=1$.

Just as in the classical case $q=1$, \eqref{8.2.1} is a resolution of the
trivial $U_q\mathfrak{g}$-module in the category $\mathcal{O}$, the full
subcategory of finitely generated weight $U_q\mathfrak{b}^+$-finite
modules.

\begin{remark}\label{8.2.9}
While constructing a BGG resolution, a function on the set of edges of the
oriented graph $\mathbb{G}$ with the values $\pm 1$, and such that the
product of its values on the edges of each square is $-1$, was used. It is
easy to establish via an argument similar to that of \ \cite{Rocha} , p.
355, 356, that a function with such properties is essentially unique. More
precisely, for any two functions $\epsilon(w_1,w)$, $\varepsilon(w_1,w)$ of
this form, one has
\begin{equation}\label{8.2.4}
\epsilon(w_1,w)=\gamma(w_1)^{-1}\varepsilon(w_1,w)\gamma(w)
\end{equation}
with a function $\gamma:W\to\{+1,-1\}$. In fact, let $W^{(k)}=\{w\in
W|\,l(w)=k\}$. For each $w\in W^{(k)}$, $k\ge 1$, fix such $w'\in
W^{(k-1)}$ that $w'\to w$ and $w=w's_\alpha$ for some simple root $\alpha$.
Define a function $\gamma(w)$ recursively:
\begin{equation}\label{8.2.5}
\gamma(e)=1,\qquad
\gamma(w)=\gamma(w')\frac{\epsilon(w',w)}{\varepsilon(w',w)}.
\end{equation}
\eqref{8.2.4} can be proved by an induction argument in $k$ which uses the
properties of $\epsilon(w_1,w)$, $\varepsilon(w_1,w)$ (Lemma
\ref{BGG-square} and Lemma 11.3 of \cite{BGG}) . It follows from
\eqref{8.2.4} that the BGG resolutions corresponding to $\epsilon(w,w')$
and $\varepsilon(w,w')$ are isomorphic in the category of complexes of
$U_q\mathfrak{g}$-modules. One can use a family of the linear maps
$\mu_k:C_k\to C_k$ with $\mu_k|_{M(w\cdot
0)}=\gamma(w)\operatorname{id}_{M(w\cdot 0)}$ as an isomorphism of
complexes.
\end{remark}

\section{\sffamily\bfseries\!\!\!\!\!\!. A GENERALIZED BGG RESOLUTION}
\label{Gen_BGG}

\indent

Some results of this Section are valid not only for
$\mathbb{S}=\{1,2,\ldots,l\}\setminus\{l_0\}$ but for any subset
$\mathbb{S}\subset\{1,2,\ldots,l\}$, and the associated Hopf algebras
$U_q\mathfrak{k}$, $U_q\mathfrak{q}^+$, a lattice $\mathscr{P}_+$ and the
generalized Verma modules $N(\mathfrak{q}^+,\lambda)$,
$\lambda\in\mathscr{P}_+$, see \cite{Rocha} . Among the results one should
mention Proposition \ref{2.1.1} and those of the present section.

We produce a resolution of the trivial module $\mathbb{C}$ in the category
$\mathcal{O}_\mathbb{S}$, that is, the full subcategory formed by
$U_q\mathfrak{k}$-finite $U_q\mathfrak{g}$-modules of the category
$\mathcal{O}$. We follow the ideas of \cite{Lep2} and \cite{Rocha} where
this problem was solved in the classical case $q=1$.

Let $\lambda\in\mathscr{P}_+$ and let $p_\lambda$ be the canonical
surjective morphism of $U_q\mathfrak{g}$-modules
$$
p_\lambda:M(\lambda)\to N(\mathfrak{q}^+,\lambda),\qquad
p_\lambda:v(\lambda)\mapsto v(\mathfrak{q}^+,\lambda).
$$

Let $\lambda,\mu\in\mathscr{P}_+$ and $f:M(\lambda)\to M(\mu)$ be a
non-zero morphism of $U_q\mathfrak{g}$-modules. If there exists a morphism
of $U_q\mathfrak{g}$-modules $\widehat{f}:N(\mathfrak{q}^+,\lambda)\to
N(\mathfrak{q}^+,\mu)$ such that $\widehat{f}p_\lambda=p_\mu f$, it is
unique and is called the standard morphism associated with $f$.

The following result is a q-analogue of \cite{Lep2} , Proposition 3.1,
and can be proved in the same way as the result of Lepowsky we refer to.

\begin{proposition}
Let $\lambda\in P$, $\mu\in\mathscr{P}_+$, and $f:M(\lambda)\to M(\mu)$ be
a morphism of Verma modules. If $p_\mu f\ne 0$, then
$\lambda\in\mathscr{P}_+$ and there exists a standard morphism of the
generalized Verma modules $\widehat{f}:N(\mathfrak{q}^+,\lambda)\to
N(\mathfrak{q}^+,\mu)$.
\end{proposition}

The subset $W^\mathbb{S}$ was used in Sec. \ref{new_item}, and the graph
$\mathbb{G}$ with $W$ being the set of vertices, was an essential point of
the previous Section. We are going to show that for one to pass to the
generalized Verma modules (as well as to the generalized BGG resolution)
one should only replace the graph by its subgraph determined by the set of
vertices $W^\mathbb{S}$.

It is known from \cite{Lep2} , p. 502, that $w\cdot\mu\in \mathscr{P}_+$
for any $\mu\in P_+$ and $w\in\,W ^\mathbb{S}$. Similarly to \eqref{8.2.1},
consider the complex of $U_q\mathfrak{g}$-modules
\begin{equation}\label{8.2.6}
0\longrightarrow C_{l\left(w_0^\mathbb{S}\right)}^\mathbb{S}
\overset{d_{l\left(w_0^\mathbb{S}\right)}}{\longrightarrow}\ldots
\overset{d_2}{\longrightarrow}C_1^\mathbb{S}\overset{d_1}{\longrightarrow}
C_0^\mathbb{S}\overset{\epsilon}{\longrightarrow}\mathbb{C}\longrightarrow
0,
\end{equation}
with $\epsilon:N(\mathfrak{q}^+,0)\to\mathbb{C}$,
$\epsilon:v(\mathfrak{q}^+,0)\mapsto 1$, being the obvious surjective
morphism,
$$
C_j^\mathbb{S}=
\bigoplus_{\{w\in\,W^\mathbb{S}\,|\:l(w)=j\}}N(\mathfrak{q}^+,w\cdot 0),
$$
and $w_0^\mathbb{S}$ being the longest element in $W^\mathbb{S}$.

Define the differentials $d_j$ just as in \eqref{8.2.3}:
\begin{equation}\label{8.2.7}
d_j|_{N(\mathfrak{q}^+,w\cdot 0)}=\bigoplus_{\{w'\in
 \,W^\mathbb{S}\,|\:w'\to\, w\}}\epsilon(w,w')\,\widetilde{i}_{w,w'},
\end{equation}
where
$$
\widetilde{i}_{w,w'}:N(\mathfrak{q}^+,w\cdot 0)\to N(\mathfrak{q}^+,w'\cdot
0),\qquad \widetilde{i}_{w,w'}:v(\mathfrak{q}^+,w\cdot 0)\mapsto p_{w'\cdot
0}i_{w,w'}v(w\cdot 0).
$$

The relation $d_j\circ d_{j+1}=0$ follows from a similar relation for
\eqref{8.2.1}, together with the following statement whose proof is similar
to that of \ \cite{Lep2} , p. 503.

\begin{lemma}
Let $w,w'\in\, W^\mathbb{S}$ and $l(w)=l(w')+1$. A non-zero morphism of
$U_q\mathfrak{g}$-modules $N(\mathfrak{q}^+,w\cdot 0)\to
N(\mathfrak{q}^+,w'\cdot 0)$ exists iff $w'\to w$.
\end{lemma}

The exactness of \eqref{8.2.6} can be proved in the same way as that of
\eqref{8.2.1} in Sec. \ref{BGG-section}, i.e., via the exactness of the
generalized BGG resolution in the classical case $q=1$ (see \cite{Lep2,
Rocha}).

In Sec. \ref{new_FODC} we showed the following fact: consider the graded
$U_q\mathfrak{g}$-module $N(\mathfrak{q}^+,\lambda)$, then its graded dual
possesses a $U_q\mathfrak{g}$-module
$\mathbb{C}[\mathfrak{p}^-]_q$-bimodule structure.

Use \eqref{8.2.6} to replace its terms by the dual graded vector spaces and
the morphisms involved therein by the adjoint linear maps to get a complex
of the $U_q\mathfrak{g}$-module $\mathbb{C}[\mathfrak{p}^-]_q$-bimodules
\begin{equation}\label{8.2.8}
0\longrightarrow\mathbb{C}\longrightarrow\left(C_0^\mathbb{S}\right)^*
\overset{d_1}{\longrightarrow}\left(C_1^\mathbb{S}\right)^*
\overset{d_2}{\longrightarrow}\ldots
\overset{d_{l\left(w_0^\mathbb{S}\right)}}{\longrightarrow}
\left(C_{l\left(w_0^\mathbb{S}\right)}^\mathbb{S}\right)^*\longrightarrow 0.
\end{equation}

\begin{proposition}
The complex of the $U_q\mathfrak{g}$-module
$\mathbb{C}[\mathfrak{p}^-]_q$-bimodules \eqref{8.2.8} is exact.
\end{proposition}

It should be noted that, while producing the complex of the
$U_q\mathfrak{g}$-module $\mathbb{C}[\mathfrak{p}^-]_q$-bimodules
\eqref{8.2.8}, a specific choice of the function $\epsilon(w,w')$ on the
set of edges of the graph $\mathbb{G}$ has been used implicitly. It follows
from the observation of the previous Section that the isomorphism class in
the category of complexes of $U_q\mathfrak{g}$-module
$\mathbb{C}[\mathfrak{p}^-]_q$-bimodules does not depend on this choice.

\section{\sffamily\bfseries\!\!\!\!\!\!. A DE RHAM COMPLEX}\label{de_Rham}

\indent

In this Section, we turn back to the assumption
$\mathbb{S}=\{1,2,\ldots,l\}\setminus\{l_0\}$, with $\alpha_{l_0}$ being a
simple root whose coefficient in \eqref{decomposition} is $1$.

We continue with our study of the differential calculus
$(\Lambda(\mathfrak{p}^-)_q,d)$ which we started in Sec.
\ref{new_diff_univ}.

We begin with the classical case $q=1$. Consider the complex of graded
$U\mathfrak{g}$-modules dual to the de Rham complex
$$
0\longrightarrow\mathbb{C}\longrightarrow
\Lambda^0(\mathfrak{p}^-)\overset{d_1}{\longrightarrow}
\Lambda^1(\mathfrak{p}^-)\overset{d_2}{\longrightarrow}\ldots
\longrightarrow\Lambda^{\dim\mathfrak{p}^-}(\mathfrak{p}^-) \longrightarrow
0.
$$
The complex is known to be isomorphic to a generalized BGG resolution. (To
prove this, use the results of Rocha-Caridi \ \cite{Rocha} , p. 364, Lemma
\ref{Kostant-lemma}, the duality of the relative Koszul complex, and the de
Rham complex of differential forms with polynomial coefficients. This
duality is well known and it was discussed in the case of differential
forms with coefficients from the algebra of formal series at the origin in
\cite{BGG}).

Now we turn to quantum analogues. Consider the universal differential
calculus $(\Lambda(\mathfrak{p}^-)_q,d)$ (see Sec. \ref{new_diff_univ}),
the linear maps $d_j=d|_{\Lambda^{j-1}(\mathfrak{p}^-)_q}$, and the de Rham
complex
\begin{equation}\label{8.2.16}
0\longrightarrow\mathbb{C}\longrightarrow
\Lambda^0(\mathfrak{p}^-)_q\overset{d_1}{\longrightarrow}
\Lambda^1(\mathfrak{p}^-)_q\overset{d_2}{\longrightarrow}\ldots
\overset{d_{\dim\mathfrak{p}^-}}
\longrightarrow\Lambda^{\dim\mathfrak{p}^-}(\mathfrak{p}^-)_q
\longrightarrow 0.
\end{equation}
Consider also the complex which is dual to \eqref{8.2.16} in the category
of graded $U_q\mathfrak{g}$-modules. Prove that it is isomorphic to
\eqref{8.2.6}.

\begin{lemma}\label{8.2.25}
In the category of weight $U_q\mathfrak{g}$-modules one has
\begin{equation}\label{8.2.17}
\left(\Lambda^k(\mathfrak{p}^-)_q\right)^*\approx \bigoplus_{\{w\in
W^\mathbb{S}|\:l(w)=k\}}N(\mathfrak{q}^+,w\cdot 0),\qquad
k=1,2,\ldots\dim\mathfrak{p}^-.
\end{equation}
\end{lemma}

{\it Proof.} Use the fact that this statement is valid for $q=1$. Choose
the weight monomial bases in $N(\mathfrak{q}^+,w\cdot 0)$,
$\left(\Lambda^k(\mathfrak{p}^-)_q\right)^*$ in the same way as in Sec.
\ref{new_item}, \ref{new_diff_univ}. The weights of vectors in these bases
do not change under quantization. This allows one to use the results of
Sec. \ref{new_diff_univ}, and the universal property of the generalized
Verma modules to prove the existence of such morphisms of
$U_q\mathfrak{g}$-modules
$$
\mathcal{J}_k(q):\bigoplus_{\{w\in\,W^\mathbb{S}\:|\:l(w)=k\}}
N(\mathfrak{q}^+,w\cdot 0)\to\left(\Lambda^k(\mathfrak{p}^-)_q\right)^*,
$$
that, firstly, their matrix elements with respect to the chosen bases are
rational functions from $\mathbb{Q}(q)$ having no poles at $q=1$, and,
secondly, the linear map $\mathcal{J}_k(q)$ is one-to-one for $q=1$. Since
the homogeneous components of the $U_q\mathfrak{g}$-modules in question are
finite dimensional, it follows that the morphism of
$U_q\mathfrak{g}$-modules $\mathcal{J}_k(q)$ is one-to-one. \hfill
$\blacksquare$

\medskip

Use \eqref{8.2.16} to replace its terms by the dual graded vector spaces
and the morphisms involved therein by the adjoint linear maps to get the
complex
\begin{equation}\label{8.2.18}
0\longrightarrow C_{l\left(w_0^\mathbb{S}\right)}^\mathbb{S}
\overset{\delta_{\dim\mathfrak{p}^-}}{\longrightarrow}\ldots
\overset{\delta_2}{\longrightarrow}C_1^\mathbb{S}
\overset{\delta_1}{\longrightarrow}C_0^\mathbb{S}
\overset{\epsilon}{\longrightarrow}\mathbb{C}\longrightarrow 0,
\end{equation}
with $\epsilon$ being the obvious surjective morphism, $\delta_j=d_j^*$,
and
$$
C_j^\mathbb{S}=
\bigoplus_{\{w\in\,W^\mathbb{S}\:|\:l(w)=j\}}N(\mathfrak{q}^+,w\cdot 0).
$$

\begin{lemma}\label{8.2.26}
The morphisms $N(\mathfrak{q}^+,w'\cdot 0)\to N(\mathfrak{q}^+,w''\cdot 0)$
of the generalized Verma modules involved in \eqref{8.2.18} are standard.
\end{lemma}

{\it Proof.} The desired morphisms of $U_q\mathfrak{g}$-modules $M(w'\cdot
0)\to M(w''\cdot 0)$ will be obtained by a duality argument. First,
introduce a $U_q\mathfrak{g}$-module algebra $\Lambda$ and an embedding of
$U_q\mathfrak{g}$-module algebras
$\Lambda(\mathfrak{p}^-)_q\hookrightarrow\Lambda$.

The Verma module $M(0)$ with zero highest weight is a graded
$U_q\mathfrak{g}^\mathrm{cop}$-module coalgebra. The dual graded
$U_q\mathfrak{g}$-module algebra $\mathbb{C}[\mathfrak{n}^-]_q$ is a
quantum analogue for the polynomial algebra on the vector space
$\mathfrak{n}^-$. Set
$$
\Lambda\overset{\mathrm{def}}{=}\mathbb{C}[\mathfrak{n}^-]_q
\otimes_{\mathbb{C}[\mathfrak{p}^-]_q}\Lambda(\mathfrak{p}^-)_q.
$$
Equip $\Lambda$ with a structure of $U_q\mathfrak{g}$-module algebra and
show that the linear map
\begin{equation}\label{trivial_embed}
\Lambda(\mathfrak{p}^-)_q\to\Lambda,\qquad \omega\mapsto1\otimes\omega,
\end{equation}
is an embedding of $U_q\mathfrak{g}$-module algebras.

We use the so called dyslectic modules over commutative algebras in braided
tensor categories, see the Appendix.

Let $\mathcal{C}^-$ be the full subcategory of weight
$U_q\mathfrak{b}^-$-finite dimensional $U_q\mathfrak{g}$-modules. This is
an Abelian braided tensor category. The algebra
$\mathbb{C}[\mathfrak{p}^-]_q$ is commutative in $\mathcal{C}^-$, and the
$\mathbb{C}[\mathfrak{p}^-]_q$-bimodules $\Lambda(\mathfrak{p}^-)_q$ and
$\mathbb{C}[\mathfrak{n}^-]_q$ are dyslectic modules over
$\mathbb{C}[\mathfrak{p}^-]_q$. By Proposition \ref{2.2.17}, the category
of dyslectic modules over $\mathbb{C}[\mathfrak{p}^-]_q$ is an Abelian
braided tensor category with the tensor product
$\otimes_{\mathbb{C}[\mathfrak{p}^-]_q}$. This allows one to equip
$\Lambda=\mathbb{C}[\mathfrak{n}^-]_q
\otimes_{\mathbb{C}[\mathfrak{p}^-]_q}\Lambda(\mathfrak{p}^-)_q$ with the
structure of $U_q\mathfrak{g}$-module algebra
$$
\Lambda\otimes\Lambda\to\Lambda\otimes_{\mathbb{C}[\mathfrak{p}^-]_q}\Lambda
\to\Lambda
$$
in a way which is standard in the quantum group theory \cite{Majid_book} ,
p. 438. In this setting, $\Lambda(\mathfrak{p}^-)_q\hookrightarrow\Lambda$
and
$$
\Lambda=\bigoplus_{j=0}^{\dim\mathfrak{p}^-}\Lambda^j,\qquad
\Lambda^j=\mathbb{C}[\mathfrak{n}^-]_q
\otimes_{\mathbb{C}[\mathfrak{p}^-]_q} \Lambda^j(\mathfrak{p}^-)_q.
$$

We shall check that in the category of $U_q\mathfrak{g}$-modules one has
\begin{equation}\label{8.2.19}
\Lambda^j\approx\bigoplus_{\{w\in W^\mathbb{S}|\:l(w)=j\}}M(w\cdot 0)^*.
\end{equation}
In fact, in the category of weight $U_q\mathfrak{h}$-modules
$$
\mathbb{C}[\mathfrak{n}^-]_q
\otimes_{\mathbb{C}[\mathfrak{p}^-]_q}N(\mathfrak{q}^+,w\cdot
0)^*\cong\mathbb{C}[\mathfrak{n}^-]_q\otimes(N(\mathfrak{q}^+,w\cdot
0)^*)_\mathrm{lowest},
$$
with $(N(\mathfrak{q}^+,w\cdot 0)^*)_\mathrm{lowest}$ being the lowest
homogeneous component of the graded vector space $N(\mathfrak{q}^+,w\cdot
0)^*$. Thus among the weights of the $U_q\mathfrak{g}$-module
$\mathbb{C}[\mathfrak{n}^-]_q
\otimes_{\mathbb{C}[\mathfrak{p}^-]_q}N(\mathfrak{q}^+,w\cdot 0)^*$ there
exists the lowest weight, and the associated weight space of which is one
dimensional. The universal property of the Verma module $M(w\cdot 0)$
implies the existence of a non-zero morphism of $U_q\mathfrak{g}$-modules
$$
M(w\cdot 0)\to\left(\mathbb{C}[\mathfrak{n}^-]_q
\otimes_{\mathbb{C}[\mathfrak{p}^-]_q}N(\mathfrak{q}^+,w\cdot
0)^*\right)^*.
$$
This morphism is unique up to a scalar multiplier. Use a duality argument
to obtain a morphism of $U_q\mathfrak{g}$-modules
\begin{equation}\label{8.2.20}
\mathbb{C}[\mathfrak{n}^-]_q\otimes_{\mathbb{C}[\mathfrak{p}^-]_q}
\Lambda^j(\mathfrak{p}^-)_q\to\bigoplus_{\{w^{-1}\in\,
W^\mathbb{S}\,|\,l(w)=j\}}M(w\cdot 0)^*.
\end{equation}

We check that it is one-to-one. Let $M(\mathfrak{k},0)$ be the Verma module
over $U_q\mathfrak{k}$ with the highest weight $0$. Use the
Poincar\'e-Birkhoff-Witt bases to introduce the isomorphisms of
$U_q\mathfrak{h}$-modules
$$
M(\mathfrak{k},0)^*\otimes\mathbb{C}[\mathfrak{p}^-]_q
\overset{\simeq}{\longrightarrow}\mathbb{C}[\mathfrak{n}^-]_q,\qquad
M(\mathfrak{k},0)^*\otimes N(\mathfrak{q}^+,w\cdot
0)^*\overset{\simeq}{\longrightarrow}M(w\cdot 0)^*
$$
and the bases of homogeneous components of the $U_q\mathfrak{g}$-modules
$$
\mathbb{C}[\mathfrak{n}^-]_q\otimes_{\mathbb{C}[\mathfrak{p}^-]_q}
\Lambda^j(\mathfrak{p}^-)_q,\qquad\bigoplus_{\{w^{-1}\in\,
W^\mathbb{S}\,|\,l(w)=j\}}M(w\cdot 0)^*.
$$
The morphism \eqref{8.2.20} maps every homogeneous component into a
homogeneous component of the same degree. It is easy to prove that, in the
above bases, the matrices of linear maps between the homogeneous
components are square matrices whose entries are in
$\mathbb{C}(q^{\frac1{s}})$ and which are regular at $q^{\frac1{s}}=1$.
Hence their invertibility at transcendental $q$ follows from their
invertibility at $q=1$. So it remains to apply the invertibility of the
morphism \eqref{8.2.20} at $q=1$. Thus we get the isomorphism
\eqref{8.2.19}.

To complete the proof of the lemma, extend the endomorphism
$d:\Lambda(\mathfrak{p}^-)_q\to\Lambda(\mathfrak{p}^-)_q$ of the
$U_q\mathfrak{g}$-module $\Lambda(\mathfrak{p}^-)_q$ to an endomorphism
$d_\mathrm{ext}:\Lambda\to\Lambda$ of the $U_q\mathfrak{g}$-module
$\Lambda$. Use a covariant first order differential calculus
$(M(-\alpha_{l_0})^*,d_{\mathfrak{p}^-})$ over the algebra
$\mathbb{C}[\mathfrak{n}^-]_q$, which can be derived just as in Sec.
\ref{new_FODC} by a duality argument from the morphism of
$U_q\mathfrak{g}$-modules
$$M(-\alpha_{l_0})\to M(0),\qquad v(-\alpha_{l_0})\mapsto F_{l_0}v(0).$$

The above isomorphism of the $U_q\mathfrak{g}$-module
$\mathbb{C}[\mathfrak{n}^-]_q$-bimodules $\Lambda^1\approx
M(-\alpha_{l_0})^*$ is determined up to a scalar multiplier, which can be
chosen in such a way that the associated operator
$d_{\mathfrak{p}^-}:\mathbb{C}[\mathfrak{n}^-]_q\to\Lambda^1$ satisfies
\begin{equation}\label{8.2.22}
d_{\mathfrak{p}^-}\varphi=d\varphi,\qquad
\varphi\in\mathbb{C}[\mathfrak{p}^-]_q.
\end{equation}

We get the desired extension $d_\mathrm{ext}$ of $d$ by setting
$$d_\mathrm{ext}(f\omega)=(d_{\mathfrak{p}^-}f)\omega+fd\omega$$
for all $f\in\mathbb{C}[\mathfrak{n}^-]_q\subset\Lambda$,
$\omega\in\Lambda(\mathfrak{p}^-)_q\subset\Lambda$. This extension is well
defined due to the universal property of tensor product in the category of
$\mathbb{C}[\mathfrak{p}^-]_q$-bimodules, together with the relation
$$
(d_{\mathfrak{p}^-}(f\varphi))\omega+f\varphi
d\omega=(d_{\mathfrak{p}^-}f)\varphi\omega+fd(\varphi\omega),
$$
with $f\in\mathbb{C}[\mathfrak{n}^-]_q$,
$\omega\in\Lambda(\mathfrak{p}^-)_q$,
$\varphi\in\mathbb{C}[\mathfrak{p}^-]_q$. This relation is an easy
consequence of \eqref{8.2.22}.

It follows from the definition that $d_\mathrm{ext}$ is a morphism of
$U_q\mathfrak{g}$-modules whose restriction to
$\Lambda(\mathfrak{p}^-)_q\subset\Lambda$ coincides with the differential
$d$. \hfill $\blacksquare$

\begin{proposition}\label{dual_BGG}
In the category of $U_q\mathfrak{g}$-modules the complex \eqref{8.2.18} is
isomorphic to the generalized Berstein-Gelfand-Gelfand resolution
\eqref{8.2.6}.
\end{proposition}

{\it Proof.} It follows from Lemmas \ref{8.2.25}, \ref{8.2.26} that the
complexes under consideration may be treated as those formed by the same
$U_q\mathfrak{g}$-modules, with the differentials having a similar form
$$
d_j|_{N(\mathfrak{q}^+,w\cdot 0)}=\bigoplus_{\{w'\in\,W^\mathbb{S}\:|\:w'\to
w\}}\epsilon(w,w')\widetilde{i}_{w,w'},
$$
where only the values of the scalar multipliers $\epsilon(w,w')$ are
different. However, the argument of Sec. \ref{BGG-section} (see Remark
\ref{8.2.9}) makes sure that a replacement of scalar multipliers reduces to
the replacement of a complex of $U_q\mathfrak{g}$-modules by an isomorphic
complex. \hfill $\blacksquare$

\begin{remark}
\begin{enumerate}
\item In a recent preprint by Heckenberger and Kolb \
    \cite{HeckenbergerKolb_last} the de Rham complex of \
    \cite{HeckenbergerKolb03} , that is very similar to ours, was
    realized as a complex dual to a BGG resolution. But the methods
    used there are quite different and the precise relation between the
    two resulting de Rham complexes still remains to be clarified.

\item With some abuse of terminology, one can say that the main results
    of the paper were obtained under the assumption of $q$ being
    transcendental. The recent bright results of Heckenberger and Kolb
    \cite{HeckenbergerKolb06} give us a hope to extend our result for
    all $q$ which are not the roots of unity.
\end{enumerate}
\end{remark}

\section*{\sffamily\bfseries ACKNOWLEDGEMENT}

\indent

The authors are grateful to Stefan Kolb for helpful discussions. Also
special thanks to the Referee for pointing out numerous faults in a
previous version of this paper.

\section*{\sffamily\bfseries APPENDIX: ALGEBRAS AND MODULES IN TENSOR
CATEGORIES}\def\thesection{A}

\indent

Consider the Abelian braided tensor category $\mathcal{C}^-$, a full
subcategory of weight $U_q\mathfrak{b}^-$-finite $U_q\mathfrak{g}$-modules.
The braiding $\check{R}_{V_1,V_2}:V_1\otimes V_2\to V_2\otimes V_1$ is
defined as usual in terms of the universal R-matrix \cite{Kassel_QG_engl} .

An algebra $F$ in the category $\mathcal{C}^-$ is said to be commutative in
this category if $m=m\check{R}_{FF}$, with $m:F\otimes F\to F$ being the
multiplication in $F$.

Consider a bimodule $E$ over a commutative algebra $F$ in the category
$\mathcal{C}^-$:
$$m_\mathrm{left}:F\otimes E\to E,\qquad m_\mathrm{right}:E\otimes F\to
E.$$ We follow \cite{LychaginPrasolov} in calling it symmetric if
$m_\mathrm{right}=m_\mathrm{left}\check{R}_{EF}$.

\begin{example}
An algebra $F$ in the category $\mathcal{C}^-$ is a bimodule over any of
its subalgebras $F'$. If $F$ is commutative in $\mathcal{C}^-$, this
bimodule is symmetric.
\end{example}

\begin{proposition}(\cite{Pareigis} , Remark 1.1)
If $E_1$, $E_2$ are symmetric bimodules over a commutative algebra $F$ in
the category $\mathcal{C}^-$, then the bimodule $E_1\otimes_F E_2$ is also
symmetric.
\end{proposition}

\begin{proposition}\label{2.2.18}(\cite{Pareigis} , Remark 1.2)
Let $E$ be a left module $m_\mathrm{left}:F\otimes E\to E$ over a
commutative algebra $F$ in $\mathcal{C}^-$. The morphism
$m_\mathrm{right}=m_\mathrm{left}\check{R}_{EF}$ equips $E$ with a
structure of symmetric bimodule over $F$ in $\mathcal{C}^-$.
\end{proposition}

This means that every left module $E$ over a commutative algebra $F$ in
$\mathcal{C}^-$ is a symmetric bimodule over $F$ in $\mathcal{C}^-$.

A bimodule $E$ over a commutative algebra $F$ in the category
$\mathcal{C}^-$ is said to be dyslectic if
$$
m_\mathrm{right}=m_\mathrm{left}\check{R}_{EF},\qquad
m_\mathrm{left}=m_\mathrm{right}\check{R}_{FE}.
$$

\begin{proposition}\label{2.2.17}(\cite{Pareigis} , Proposition 2.4,
Theorem 2.5) Let $E_1$, $E_2$ be dyslectic bimodules over a commutative
algebra $F$ in the category $\mathcal{C}^-$ and $j$ be a canonical
surjective morphism $E_1\otimes E_2\to E_1\otimes_F E_2$. Then
\begin{enumerate}
\item $E_1\otimes_F E_2$ is a dyslectic bimodule over $F$;

\item there exists a unique morphism
    $\overline{R}_{E_1E_2}:E_1\otimes_F E_2\to E_2\otimes_F E_1$ in the
    category $\mathcal{C}^-$ such that
    $\overline{R}_{E_1E_2}j=j\check{R}_{E_1E_2}$.
\end{enumerate}
\end{proposition}

It follows in a natural way from the above that the category of dyslectic
bimodules over a commutative algebra in the category $\mathcal{C}^-$ is an
Abelian braided tensor category.

Proofs of Propositions \ref{2.2.17} and \ref{2.2.18} were obtained by
Pareigis who used the methods of category theory in a very wide generality.
Note that the modules over commutative algebras in the braided Abelian
tensor categories are used in algebraic K-theory \cite{LychaginPrasolov}
and arise naturally in the conformal quantum field theory
\cite{KirillovOstrik, Fuchs} .

\end{document}